\documentclass[10pt]{article}

\title{Explicit formula of a supersingular polynomial for rank-$2$ Drinfeld modules and applications}

\author{
Takehiro HASEGAWA\footnote{The author was supported by JSPS KAKENHI Grant Number 15K17508} \footnote{thasegawa3141592@yahoo.co.jp} \\ % (up date 12/May/2017) \\
Shiga University, Otsu, Shiga 520-0862, Japan \\
}

\usepackage{makeidx}
\usepackage{amsmath}
\usepackage{amssymb}
\usepackage{amsthm}
\usepackage{enumerate}
\usepackage[all]{xy}

\makeindex

\newtheorem{thm}{Theorem}[section]

\newtheorem{prop}[thm]{Proposition}
\newtheorem{cor}[thm]{Corollary}
\newtheorem*{maintheorem*}{Main theorem}
\newtheorem*{cor*}{Corollary}
\newtheorem*{Keylemma*}{Keylemma}
\newtheorem*{fact*}{Fact}

\theoremstyle{definition}
\newtheorem*{abstract*}{Abstract}
\newtheorem*{Def*}{Definition}
\newtheorem*{rem*}{Remark}
\newtheorem*{claim*}{Claim}
\newtheorem{ex}[thm]{Example}
\newtheorem*{acknowledgements*}{Acknowledgements}

\newcommand{\N}{{\mathbb{N}}}
\newcommand{\Z}{{\mathbb{Z}}}

\newcommand{\F}{{\mathbb{F}}}

\begin{document}

\maketitle

%%%%%%%%%%%%%
\begin{abstract}
Rank-$2$ Drinfeld modules are a function-field analogue of elliptic curves, and the purpose of this paper is to investigate similarities and differences between rank-$2$ Drinfeld modules and elliptic curves in terms of supersingularity. 
Specifically, we provide an explicit formula of a supersingular polynomial for rank-$2$ Drinfeld modules and prove several basic properties. 
As an application, we give a numerical example of an asymptotically optimal tower of Drinfeld modular curves. 
\end{abstract}
%%%%%%%%%%%



%%%%%%%%%%%%%%%%%%
\section{Introduction}\label{sectintro}
%%%%%%%%%%%%%%%%%%

Throughout this paper, we use the following terminology. 
We let $\Z$ denote the set of integers, and $\N$ denote the set of nonnegative integers. 
That is, 
$$
\N := \{0, 1, 2,3,  \ldots\}. 
$$
The cardinality of a finite set $S$ is denoted by $|S|$. 
The notation $\lceil x \rceil$ indicates the ceiling of $x$ and denotes the smallest integer greater than or equal to a real number $x$, namely, $\lceil x \rceil := \min \{n \in \Z \ | \ x \leq n \}$.  
In this paper, $q$ will be some power of a prime number, and $\F_{q}$ will denote the finite field with $q$ elements. 
The algebraic closure of a field $L$ is denoted $\bar{L}$. 
% A fixed algebraic closure of a field $L$ is denoted $\bar{L}$. 

\medskip

Throughout this paper, we adopt the notation used in the following references: Goss \cite{Gos96} and Thakur \cite{Tha04} for rank-$2$ Drinfeld modules, Gekeler \cite{Gek86, Gek01, Gek04} for Drinfeld modular curves, Silverman \cite{Sil09} and Husem\"{o}ller \cite{Hus04} for elliptic curves, and Stichtenoth \cite{Sil09} for function fields. 

\bigskip

% \medskip

It is known that rank-$2$ Drinfeld modules are a function-field analogue of elliptic curves. 
This relation was first discovered by Drinfeld \cite{Dri74, Dri77} and has been studied since then by many researchers (see, for example, \cite{Gek86}, \cite{Gos96}, \cite{Tha04}). 
It is thus natural to investigate similarities and differences between rank-$2$ Drinfeld modules and elliptic curves. 
This paper does so in terms of supersingularity. 
% This paper is part of such investigations. 
% In this paper, we prove 

\bigskip

% \medskip

Before we present the main theorem and its corollary, we recall the analogous results for the elliptic-curve case. % , and recent results . 
% Let $\C$ denote the set of complex numbers, and set $\lambda \in \C \setminus \{0,1 \}$. 
Let $p \geq 3$ be a prime number. %, and let $\bar{\F}_{p}$ denote an algebraic closure of $\F_{p}$. 
% Let $\F_{p}:= \Z/(p)$ denote the finite field, and let $\bar{\F}_{p}$ denote an algebraic closure of $\F_{p}$. 
It is well-known that every elliptic curve is isomorphic (over $\bar{\F}_{p}$) to an elliptic curve in Legendre  form 
$$
E_{\lambda} \colon y^{2}= x(x-1)(x- \lambda), 
$$
where $\lambda$ is an element in $\bar{\F}_{p}$ with $\lambda \neq 0, 1$ (see, for example, Proposition 1.7 of Chapter III in \cite{Sil09}). 
Let $m \geq 1$ be an integer, and let $E_{\lambda}(\bar{\F}_{p})[m]$ denote the $m$-torsion subgroup of $E_{\lambda}$. % , that is, the set of points of order $m$ in $E_{\lambda}$. 
The elliptic curve $E_{\lambda}$ is called supersingular when $E_{\lambda}(\bar{\F}_{p})[p]=0$. 

From here, $\lambda$ is regarded as an indeterminate element. 
We set the Deuring polynomial
$$
H_{p}(\lambda):= \sum_{i=0}^{(p-1)/2} \binom{(p-1)/2}{i}^{2} \lambda^{i} \in \F_{p}[\lambda], 
$$
which was first defined by Deuring \cite{Deu41}. 
Then, we know the following facts: 

\medskip

(E1) \ 
$E_{\lambda}$ is supersingular if and only if $H_{p}(\lambda)=0$ (e.g., see Theorem 4.1 (b) of Chapter V in \cite{Sil09}). 

\medskip

(E2) \ 
$H_{p}(\lambda)$ is separable (e.g., see Theorem 4.1 (c) of Chapter V in \cite{Sil09}). 

\medskip

(E3) \ 
If $\lambda \in \bar{\F}_{p}$ is a root of $H_{p}(\lambda^{4})=0$, then $\lambda \in \F_{p^{2}}$ (see Theorem of Appendix in \cite{GSR03}). 
% In particular, the roots of $H_{p}(\lambda)$ are in $\F_{p^{2}}$. 

\medskip

(E3$_\text{weak}$) \ 
If $\lambda \in \bar{\F}_{p}$ is a root of $H_{p}(\lambda)=0$, then $\lambda \in \F_{p^{2}}$. 

\medskip

(E4) \ 
The sequence of elliptic modular curves $X_{0}(2^{n})/\F_{p^{2}}$ ($n \geq 2$) is asymptotically optimal, which is expressed as 
$$
\frac{N(X_{0}(2^{n})/\F_{p^{2}})}{g(X_{0}(2^{n}))} \to p-1 \qquad (n \to \infty),  
$$
where $N(X_{0}(2^{n})/\F_{p^{2}})$ denotes the number of rational points of $X_{0}(2^{n})/\F_{p^{2}}$, and $g(X_{0}(2^{n}))$ denotes the genus of $X_{0}(2^{n})$ (for a more elementary proof, see Theorem 5.2 in \cite{GSR03}, and moreover, for a proof of a more general case, see Theorem 3.2 in \cite{Gek04} and Theorem 4.1.52 in \cite{TV91}). 

\medskip

Note that although a more general case of (E4) is proven, the proof in \cite{GSR03} explicitly describes the set of rational points that split completely in the sequence, in order to apply coding theory. 

\bigskip

We know the following function-field analogues of (E1), (E2), (E3) and (E3$_\text{weak}$). 
Specifically, Gekeler has proven analogues of (E1) and (E3$_\text{weak}$) for rank-$2$ Drinfeld modules defined by the general type (\ref{drinfeldgeneral}) (see Satz (5.3) in \cite{Gek83}, Proposition 4.2 in \cite{Gek91}), and has proven analogues of (E2) and (E3) for a rank-$2$ Drinfeld module defined by the specific type (\ref{gekelernormal}) (see Lemmas (5.6) and (5.7) in \cite{Gek83}), where the types (\ref{drinfeldgeneral}) and (\ref{gekelernormal}) are defined later. 
Notice that Proposition 4.2 in \cite{Gek91} discusses Drinfeld modules of arbitrary rank. 
El-Guindy and Papanikolas have computed a formula that corresponds to $H_{p}(\lambda)$ for Drinfeld modules of arbitrary rank (see Theorem 8.1, Corollary 8.2 in \cite{EP13}), but it is not explicit in our sense. 
Bassa and Beelen have defined polynomials by recursion, and have proven analogues of (E1), (E2) and (E3$_\text{weak}$) for these polynomials (see Corollary 20, Proposition 13, Theorem 18 in \cite{BB}). 
El-Guindy has computed an explicit formula for the polynomials of Bassa and Beelen (see Theorem 3.2 (ii) in \cite{ElG14}). 
From the viewpoint of the theory of Drinfeld modules, the results of Bassa and Beelen can be regarded as analogues of (E1), (E2) and (E3$_\text{weak}$) for a rank-$2$ Drinfeld module defined by the specific type (\ref{drinfeldnormal}) (see the remark in Section \ref{sectss} of this paper), where the type (\ref{drinfeldnormal}) is defined later. 
Similarly, the result of El-Guindy can be regarded as an explicit formula that corresponds to $H_{p}(\lambda)$ for a rank-$2$ Drinfeld module defined by (\ref{drinfeldnormal}) (see the remark in Section \ref{sectss} of this paper). 

We know the following function-field analogues of (E4). 
Bassa and Beelen have studied an analogue of (E4) using analogues of (E1), (E2) and (E3$_\text{weak}$) (see Corollary 17 in \cite{BB}). 
An analogue of a more general case of (E4) has also been studied (see Theorem 2.16 in \cite{Gek04}, Theorem 4.2.38 in \cite{TV91}). 

\medskip

% To my knowledge, 
An explicit formula that corresponds to $H_{p}(\lambda)$ for rank-$2$ Drinfeld modules defined by the general type (\ref{drinfeldgeneral}) has not been before this paper. % (cf. Section 5 in \cite{Gek83}, Corollary 8.2 in \cite{EP13}). 
In the first half of this paper, we provide it, with equivalence proven as Main theorem (1). 
% This is Main theorem (1). 
Using Main theorem (1), we prove analogues of (E2) and (E3) for a rank-$2$ Drinfeld module defined by (\ref{drinfeldnormal}) (as Main theorem (2)). 
In the final half of this paper, we prove an analogue of (E4) as a corollary by applying the theorem. 

Although Main theorem (2) is similar to Lemmas (5.6) and (5.7) in \cite{Gek83}, we prove it here because we need it in the proof of the corollary. 
Note that our corollary has already been studied by Bassa and Beelen (see Corollary 17 in \cite{BB}). 
However,  their proof becomes complicated in showing that the completely splitting points are rational (see Pages 12 and 13), because they do not use a supersingular polynomial in the Gekeler's sense. 
In this paper, the proof of this is simple, because we use a supersingular polynomial. 
Note that our corollary is a special case of Theorem 2.16 in \cite{Gek04} and Theorem 4.2.38 in \cite{TV91}. 
However, our result has interesting applications within coding theory (see \cite{Sti06}, Chapters 7 and 8 in \cite{Sti09}, Parts 3 and 4 in \cite{TV91}, Chapters 3 and 4 in \cite{TVD07}), because our proof is more elementary, and explicitly describes the set of rational points that split completely. 

\bigskip

The low-genus curve cases of (E1)--(E4) were studied in \cite{Has10, Has13, Has15}. 
There, we defined explicit polynomials that correspond to $H_{p}(\lambda)$ and showed that the sequences of elliptic modular curves $X_{0}(3^{n})$, $X_{0}(4^{n})$, $X_{0}(3 \cdot 2^{n})$, and $X_{0}(2 \cdot 3^{n})$ ($n \geq 2$) are asymptotically optimal. %, respectively. 

\bigskip

% \medskip

Thakur introduced two distinct hypergeometric functions for function fields (see \cite{Tha95, Tha00}, Subsection 6.5 in \cite{Tha04}). 
In this paper, we suggest the possibility of another hypergeometric function. % (see the last remark of Section \ref{sectss}). 
As background material, we recall a relation between the polynomial $H_{p}(\lambda)$, a hypergeometric function, and a period of an elliptic curve. 
It is known that a (real) period 
$$
\omega_{2}(\lambda) := \int_{1}^{\infty} \frac{dx}{\sqrt{x(x-1)(x- \lambda)}} 
$$
of the elliptic curve $E_{\lambda}$ is equal to the product of the Gauss hypergeometric function 
$$
F \left( {1}/{2}, {1}/{2}, 1; \lambda \right) := \sum_{n \geq 0} \frac{(1/2)_{n} (1/2)_{n}}{(1)_{n}} \cdot \frac{\lambda^{n}}{n!} 
$$
and the circular constant $\pi$ (the ratio of circumference to diameter), that is, 
\begin{align}\label{omegaphiF}
\omega_{2}(\lambda)= \pi \cdot F \left( {1}/{2}, {1}/{2}, 1; \lambda \right)
\end{align}
(see, for example, Theorem (6.1) in \cite{Hus04}). 
For comparison, $H_{p}(\lambda)$ is the truncated hypergeometric function 
$$
F^{\text{tr}} \left( {1}/{2}, {1}/{2}, 1; \lambda \right) := \sum_{n=0}^{{(p-1)}/{2}} \frac{(1/2)_{n} (1/2)_{n}}{(1)_{n}} \cdot \frac{\lambda^{n}}{n!}, 
$$
namely, 
\begin{align}\label{HFtr}
H_{p}(\lambda) \equiv F^{\text{tr}} \left( {1}/{2}, {1}/{2}, 1; \lambda \right) \pmod p 
\end{align}
(see Page 261 in \cite{Hus04}). 
Hence, by applying the equalities (\ref{omegaphiF}) and (\ref{HFtr}), we can regard $H_{p}(\lambda)$ as the product of a suitable period $\omega_{2}^{\text{tr}}(\lambda)$ and a suitable constant $1/\pi^{\text{tr}}$, that is, 
$$
H_{p}(\lambda) \equiv \frac{1}{\pi^{\text{tr}}} \cdot \omega_{2}^{\text{tr}}(\lambda) \pmod p. 
$$
Surprisingly, this phenomenon also occurs in our Drinfeld-module case (see the last remark of Section \ref{sectss}). 

\bigskip

% \medskip

For the main theorem, we introduce a rank-$2$ Drinfeld module and a partition of a subset of $\N$. 

% Let $q$ be a power of a prime number, and 
% Let $\F_{q}$ be the finite field with $q$ elements. 
Let $A:= \F_{q}[T]$ denote a polynomial ring, and let $\mathfrak{p}$ denote its nonzero prime ideal. 
% Since $A$ is a principal ideal domain, 
Then, there is a monic irreducible polynomial $p(T) \in A$ such that $\mathfrak{p}= (p(T))$. 
Throughout this paper, we always suppose that $p(T) \neq T$. 
Set $\F_{\mathfrak{p}} :=A/\mathfrak{p}$ and $d:= \deg_{T} p(T)$. 
Let $\F_{\mathfrak{p}}^{(2)}$ denote the quadratic extension of $\F_{\mathfrak{p}}$. 
Let $\alpha$ be any root of $p(T)$, and fix this root. 
Note that $\alpha \neq 0$. % since $p(T) \neq T$.  
%  let $\mathfrak{p}= (p(T))$ denote a nonzero prime ideal of $A$, where  $p(T):= \sum_{i=0}^{d} \mu_{i} T^{i} \in A$ is a monic irreducible polynomial. 
We see that $A/\mathfrak{p}= \F_{q^{d}}= \F_{q}(\alpha)$. % from a theory of fields. 
Observe that $\alpha^{q}, \ldots, \alpha^{q^{d-1}}$ are the other roots of $p(T)$, and that $\alpha^{q^{d}}= \alpha$. 

% Throughout this paper, we fix
% Note that $\alpha^{q^{i}}$ is also a root of $p(T)$, and that $\alpha^{q^{d}}= \alpha$. 
% In fact, $p(\alpha^{q^{i}})= p(\alpha)^{q^{i}}=0$. 

Let $K:= \F_{q}(T)$ denote the quotient field of $A$, and let $C_{\infty}$ denote the completion of an algebraic closure of  the completion of $K$ at the infinite place (see Subsection 4.1 in \cite{Gos96}). 
% Let $K_{\infty}= K((1/T))$ be the completion of $K$ at the infinite place, and let $C_{\infty}$ denote the completion of an algebraic closure of $K_{\infty}$.  

\medskip

Let $L$ be an extension of either $K$ or $\F_{\mathfrak{p}}$, and let $\iota \colon A \to L$ denote an $\F_{q}$-algebra homomorphism. 
Notice that if $L$ is an extension of $K$ (resp. $\F_{\mathfrak{p}}$), then $\iota(T)=T$ (resp. $\iota(T)= \alpha$). 
Let 
$$
\tau: L \to L, \qquad \tau(l)= l^{q}
$$
denote a Frobenius endomorphism, and let $L \{ \tau \}$ denote a polynomial ring in $\tau$ under addition and composition, that is, $\tau l= l^{q} \tau$ for any $l \in L$ (see Section 1 of Chapter I in \cite{Gek86} or Subsection 1.1 in \cite{Gos96}). 

A rank-$2$ Drinfeld module over $L$ is an $\F_{q}$-algebra homomorphism 
$$
\phi: A \to L \{ \tau \}
$$
such that 
\begin{align}\label{drinfeldgeneral}
% \phi_{T}:= \phi(T):= \iota(T) \tau^{0}+ A \tau+ B \tau^{2}, \quad A,B \in L, \quad B \neq 0. 
\phi_{T} := \phi(T) = \iota(T) + A_{1} \tau + A_{2} \tau^{2}
\end{align}
($A_{1}, A_{2} \in L, A_{2} \neq 0$) (see Definition 1.3 of Chapter I in \cite{Gek86}, Definition 4.4.2 in \cite{Gos96}, Definition 2.2.1 in \cite{Tha04}). 
% In this paper, we deal with a rank $2$ Drinfeld module. 
The $j$-invariant of $\phi$ is defined by $j(\phi):= A_{1}^{q+1}/A_{2}$ (see Example 3.6 of Chapter V in \cite{Gek86}, Subsection 2.7 in \cite{Maz97}, Subsection 6.1 in \cite{Tha04}). 
Recall that for any $a \in A$, the constant term of $\phi_{a}:= \phi(a)$ is $\iota(a)$, and that the degree of $\phi_{a}$ in $\tau$ is $2 \deg_{T}(a)$. 
% It is known that two Drinfeld modules $\phi$ and $\phi^{\prime}$ are isomorphic over an algebraic closure $\bar{L}$ of $L$ if and only if $j(\phi)= j(\phi^{\prime})$. 
Let 
$$
\text{Ker}(\phi_{a}):= \{x \in \bar{L} \ | \ \phi_{a}(x)=0 \}
$$
denote the $a$-torsion points of $\phi$, which is a subspace of $\bar{L}$. 

\bigskip

% \medskip

Here, we discuss a normal form for rank-$2$ Drinfeld modules (see Section 2 in \cite{ElG14}). 
Suppose that $L= C_{\infty}$. 
In this case, it is known that $\dim_{\F_{q}} (\text{Ker}(\phi_{T}))= 2$ (see Proposition 1.6 of Chapter I in \cite{Gek86}). 
% Let $v$ denote the unique valuation on $C_{\infty}$ defined by $v(1/T)= 1$. 
% It is known that $v$ has a unique extension to $C_{\infty}$. 
Let $\delta \in \F_{q^{2}} \setminus \F_{q}$ be an element such that $\delta^{q}= - \delta$. 
Observe that $\F_{q} \delta$ is a subspace of $C_{\infty}$, and $\dim_{\F_{q}} \F_{q} \delta= 1$. 
Next, we set 
\begin{align*}
\mathcal{F}_{\delta} &:= \Big\{ \phi \ \Big| \ \text{$\phi$ is a rank-$2$ Drinfeld module over $C_{\infty}$, and $\F_{q} \delta \subseteq \text{Ker}(\phi_{T})$} \Big\}. %, \quad \text{and} \\
% \mathcal{F}^{\star}_{\delta} &:= \{ \phi \in \mathcal{F} \ | \ \text{$v(j)< -q$ and $v(A)= v(B)$} \}.  
\end{align*}
Assume that $\phi$ is any Drinfeld module in the set $\mathcal{F}_{\delta}$ defined by $\phi_{T}= \iota(T)+ A_{1} \tau+ A_{2} \tau^{2}$. 
Then, we have the relation $\lambda:= A_{2}= A_{1}- \iota(T)$.  
Conversely, any rank-$2$ Drinfeld module $\phi$ defined by 
\begin{align}\label{drinfeldnormal}
\phi_{T} = \iota(T)+ (\iota(T)+ \lambda) \tau+ \lambda \tau^{2}
\end{align}
is in the set $\mathcal{F}_{\delta}$. 
In fact, for any $x \in \F_{q}$, we have 
\begin{align*}
\phi_{T}(x \delta) &= \iota(T)(x \delta)+ ( \iota(T)+ \lambda) (x \delta)^{q}+ \lambda (x \delta)^{q^{2}} \\
&= \iota(T) x \delta- ( \iota(T)+ \lambda)x \delta+ \lambda x \delta=0.  
\end{align*}
% It is known that any Drinfeld module in the set $\mathcal{F}_{\delta}$ is given by (\ref{drinfeldnormal}) (see Section 2 in \cite{ElG14}). 
Hence, we can regard a Drinfeld module defined by (\ref{drinfeldnormal}) as being in a normal form in the above sense.  % in the above sense. 
% It follows from $j(\phi)= j(\psi)= (T+ \lambda)^{q+1}/\lambda$ that the special Drinfeld module $\phi$ given by  $$ \phi_{T}:= T+ (T+ \lambda) \tau+ \lambda \tau^{2} $$ and $\psi$ is isomorphic. 
In this paper, we sometimes restrict Drinfeld modules to those defined by (\ref{drinfeldnormal}). 

In a 1983 paper, Gekeler regarded a Drinfeld module defined by 
\begin{align}\label{gekelernormal}
\phi_{T} = \iota(T)+ \lambda \tau+ \tau^{2}
\end{align}
as being in a normal form, and using this, proved results similar to Main theorem (2) (see Section 5 in \cite{Gek83}). 
Notice that Corollary to Main theorem (2) can also be shown by using his results, though an equation used in its proof is different from our equation (\ref{BBequation}), which is defined in Section \ref{secttower}. 

\bigskip

% \medskip

For $L$ an extension of $\F_{\mathfrak{p}}$, we write 
$$
\phi_{\mathfrak{p}}:= \phi_{p(T)}= \sum_{i=0}^{2d} g_{i} \tau^{i} \qquad \text{and} \qquad H_{\mathfrak{p}}^{(d)}(\phi):= g_{d}. 
$$
Then, we know that 
$$
g_{0}= g_{1}= \cdots= g_{d-1}= 0
$$
(see Section 5 in \cite{Gek83}, Section 11 in \cite{Gek88}). 
This fact is often used in this paper. 
Recall that $\phi_{\mathfrak{p}}(0)= 0$, since $\phi_{\mathfrak{p}}(x)= \sum_{i=0}^{2d} g_{i} x^{q^{i}}$. 
A rank-$2$ Drinfeld module $\phi$ over $L$ is called supersingular at $\mathfrak{p}$ when $\text{Ker}(\phi_{\mathfrak{p}})= \{ 0 \}$ (see Remark 2.4 of Chapter VIII in \cite{Gek86}, Definition 4.12.16 in \cite{Gos96}). 
With this, $\phi$ is supersingular at $\mathfrak{p}$ if and only if $g_{d}= H_{\mathfrak{p}}^{(d)}(\phi)=0$ (see Satz (5.3) in \cite{Gek83}). 

Assume that $A_{1}= \alpha+ \lambda$ and $A_{2}= \lambda$, that is, that a Drinfeld module $\phi$ in the form (\ref{drinfeldnormal}) can be defined. 
Then, the coefficients $g_{d}, g_{d+1}, \ldots, g_{2d}$ are polynomials in $\lambda$ over $L$. 
We set 
$$
H_{\mathfrak{p}}^{(d)}(\lambda):= H_{\mathfrak{p}}^{(d)}(\phi). 
$$
% It is known that the degree of $H_{\mathfrak{p}}^{(d)}(\lambda)$ in $\lambda$ is $(q^{d}-1)/(q-1)$ (cf. Lemma (5.7) in \cite{Gek83}).  

\bigskip

% \medskip

Next, we introduce a partition of a subset of $\N$. 
For a positive integer $d>0$, we write 
$$
\N_{<d}:= \{0, 1, \ldots, d-1\}. 
$$
For a finite subset $S$ of $\N$ and a positive integer $j>0$, we let $S+j := \{i+ j \ | \ i \in S \}$. 
Notice that $\emptyset+j= \emptyset$. 

A partition of $\N_{<d}$ is a collection $\{ S_{1}, S_{2}, S_{3} \}$ of subsets of $\N_{<d}$ such that 
$$
S_{1} \cap S_{2}= S_{2} \cap S_{3}= S_{3} \cap S_{1}= \emptyset \quad \text{and} \quad S_{1} \cup S_{2} \cup S_{3}= \N_{<d}. 
$$
For an integer $d$, we define 
\begin{align*}
& P(d) = P_{2}(d) \\
& := 
\begin{cases}
\ \Big\{(S_{1}, S_{2}) \ \Big| \ \text{$\{S_{1} S_{2}, S_{2}+1 \}$ forms a partition of $\N_{<d}$} \Big\} & \text{if $d>0$}; \\
\ \{ \emptyset \} & \text{if $d=0$}; \\
\ \emptyset & \text{if $d<0$}. 
\end{cases}
\end{align*}
Note that $d= |S_{1}|+ 2 |S_{2}|$ and thus $|S_{1}|+ |S_{2}|= d- |S_{2}|$ for $(S_{1}, S_{2}) \in P(d)$. 
It is known that 
\begin{align}\label{recursionpd}
|P(d)| = 
\begin{cases}
\ |P(d-1)|+ |P(d-2)| & \text{if $d>0$}; \\
\ 1 & \text{if $d=0$}; \\
\ 0 & \text{if $d<0$} 
\end{cases}
\end{align}
(see Lemma 2.1 (iii) in \cite{EP13}). 

For a nonnegative integer $n$ and a finite subset $S$ of $\N$, set 
\begin{align*}
[n] &:= 
\begin{cases}
\ T^{q^{n}}- T & \text{if $L$ is an extension of $K$}; \\
\ \alpha^{q^{n}}- \alpha & \text{if $L$ is an extension of $\F_{\mathfrak{p}}$}, 
\end{cases} \\
L(S) &:= 
\begin{cases}
\displaystyle \ (-1)^{|S|} \prod_{i \in S} [i] & \text{if $S \neq \emptyset$}; \\
\ 1 & \text{if $S= \emptyset$}, 
\end{cases} \\ % \quad \text{and} \quad 
w(S) &:= 
\begin{cases}
\ \displaystyle \sum_{i \in S} q^{i} & \text{if $S \neq \emptyset$}; \\
\ 0 & \text{if $S= \emptyset$}. 
\end{cases}
\end{align*}

\bigskip

% \medskip

The following is our main theorem and a corollary of that theorem. 

%%%%%%%%%%%%%%%%%%
\begin{maintheorem*} %(???)
% Let $d \geq 0$ be an integer. 
Let $\mathfrak{p}= (p(T))$ denote a nonzero prime ideal of $A= \F_{q}[T]$ such that $p(T) \neq T$, and let $\alpha$ be a root of $p(T)$. 
Set $d= \deg_{T} p(T)$. 
Further, assume that $L$ is an extension of $\F_{\mathfrak{p}}$ and let $\phi$ be any rank-$2$ Drinfeld module over $L$ defined by $\phi_{T}= \alpha+ A_{1} \tau+ A_{2} \tau^{2}$. 
Then, the following hold. 

\medskip

(1) \ 
The coefficient $H_{\mathfrak{p}}^{(d)}(\phi)$ is given by 
$$
H_{\mathfrak{p}}^{(d)}(\phi) = \sum_{(S_{1}, S_{2}) \in P(d)} L(S_{2}+ 1) A_{1}^{w(S_{1})} A_{2}^{w(S_{2})}. 
$$

\medskip

(2) \ 
When $A_{1}= \alpha+ \lambda$ and $A_{2}= \lambda$, the polynomial $H_{\mathfrak{p}}^{(d)}(\lambda)$ is separable, and its degree in $\lambda$ is $\deg_{\lambda} H_{\mathfrak{p}}^{(d)}(\lambda)= (q^{d}-1)/(q-1)$. 
Moreover, if $H_{\mathfrak{p}}^{(d)}(\lambda^{q+1})=0$, then $\lambda \in \F_{\mathfrak{p}}^{(2)}$. 
\end{maintheorem*}
%%%%%%%%%%%%%%%%

\bigskip

%%%%%%%%%%
\begin{cor*}
The sequence of Drinfeld modular curves $X_{0}(T^{n})/\F_{\mathfrak{p}}^{(2)}$ ($n \geq 2$) is asymptotically optimal, that is, 
$$
\frac{N(X_{0}(T^{n})/\F_{\mathfrak{p}}^{(2)})}{g(X_{0}(T^{n}))} \to q^{d}-1 \qquad (n \to \infty), 
$$
where $N(X_{0}(T^{n})/\F_{\mathfrak{p}}^{(2)})$ denotes the number of rational points of $X_{0}(T^{n})/\F_{\mathfrak{p}}^{(2)}$, and $g(X_{0}(T^{n}))$ denotes the genus of $X_{0}(T^{n})$. 
\end{cor*}
%%%%%%%%

Note that we prove the corollary in terms of function fields, rather than curves. 
The motivation for the corollary comes from coding theory. 
For applications to the theory, it is essential that the proof is elementary and explicit. 
Our proof is elementary and explicit. % \medskip
It is well-known that computer science uses fields of characteristic $p=2$. % in many applications. 
% That is, the characteristic of the base field is $p=2$. 
Then, the base field of (E4) is $\F_{4}$, which is small. 
However, with the corollary, we can choose a base field of characteristic $p=2$ large enough. % because there exists an irreducible polynomial over $\F_{2}$ of arbitrary degree. 

\bigskip

% \medskip

The organization of this paper is as follows. 
In Section \ref{sectss}, we prove Main theorem (1) (Proposition \ref{mainthm1} (a)). 
In the proof of Main theorem (1), a polynomial identity (Keylemma) plays a key role. 
In the last part of Section \ref{sectss}, we suggest the possibility of a hypergeometric function for function fields (the last remark of Section \ref{sectss}). 
In Section \ref{applications}, we prove Main theorem (2) (Proposition \ref{mainthm23} (b), (c), (d)), which is a function-field analogue of the polynomial $H_{p}(\lambda)$. 
In Section \ref{secttower}, we prove a corollary to Main theorem (2) (as Proposition \ref{maincor}). 
In the proof of the corollary, another polynomial identity (Proposition \ref{polyidentity} (b)) plays a key role. 
In the last part of Section \ref{secttower}, we present a structure for the sequence of Drinfeld modular curves $X_{0}(T^{n})$, which is due to Sections 2 and 3 in \cite{BB}. 

% Perids, $T \log_{\phi}(\delta), T \log_{\phi}(\zeta)$

%%%%%%%%%%%%%%%%%%%%%%%%%%%%%%%%%%%%%%%%%%%%%%%
\section{An explicit formula for a supersingular polynomial}\label{sectss}  
%%%%%%%%%%%%%%%%%%%%%%%%%%%%%%%%%%%%%%%%%%%%%%%

In this section, we prove Main theorem (1) (Proposition \ref{mainthm1} (a)), which was introduced in Section \ref{sectintro}. 
%In course of the proof, a polynomial identity (Keylemma) palys very important. 
In the course of the proof, a polynomial identity (Keylemma) plays a key role. 
The proof of Main theorem (1) relies on combining our Keylemma with results of El-Guindy and Papanikolas \cite{EP13}. 

\bigskip

% \medskip

We recall the setup introduced in Section \ref{sectintro}. 
Let $\mathfrak{p}= (p(T))$ denote a nonzero prime ideal of $A= \F_{q}[T]$, where $p(T) = \sum_{i=0}^{d} \mu_{i} T^{i} \in A$ is a monic irreducible polynomial of degree $d$. 
We call $\mathfrak{p}= (p(T))$ an ideal of degree $d$. 
% Throughout this paper, we 
Let $\alpha$ be a root of $p(T)$. %, and fix this root. 
Suppose that $p(T) \neq T$. 
Note that $\alpha \neq 0$. 
Let $L$ be an extension of either $K$ or $\F_{\mathfrak{p}}$, and let $\phi$ be a Drinfeld module over $L$ defined by 
$$
\phi_{T} = \iota(T)+ A_{1} \tau+ A_{2} \tau^{2}
$$
($A_{1}, A_{2} \in L, A_{2} \neq 0$). 
% Recall that $\phi$ is a homomorphism. 

\bigskip

% \medskip

In order to understand the structure of the proof of Main theorem (1), we consider a relation between the coefficient $H_{\mathfrak{p}}^{(d)}(\phi)$ and the set $P(d)$. 
These were introduced in Section \ref{sectintro}. 
%, that is, the terms of $g_{d}$ correspond to the elements of $P(d)$. 
In the following examples, we calculate the coefficients $H_{\mathfrak{p}}^{(d)}(\phi)$ by using the rule $\tau l= l^{q} \tau$ ($l \in L$), and compute the sets $P(d)$ by using the recursion (\ref{recursionpd}) given in Section \ref{sectintro}. % ($d= 1, 2, 3, 4, 5)$. 

Notice that for $d=2,3,4,5$, our coefficients $H_{\mathfrak{p}}^{(d)}(\phi)$ and the supersingular polynomials $P_{d}(j)$ in Examples (2.2) of \cite{Cor99} coincide (cf. Proposition (6.2) (ii) (Deligne's congruence) in \cite{Cor99}). 
% Note that the item (5) is used after the proof of Keylemma. 
% Note that $H_{\mathfrak{p}}^{(5)}(\phi)$ is used after the proof of Keylemma. 
% Let $\N$ denote the set of nonnegative integers, that is, $\N := \{0, 1, 2,3,  \ldots\}$, and let $\Z$ denote the set of integers. 

% \medskip

%%%%%%%%
\begin{ex}\label{exampless} % (???) %[cf. Examples (2.2) in \cite{Cor99}]
Let $L$ be an extension of $\F_{\mathfrak{p}}$. 
Note that $[n] = \alpha^{q^{n}}- \alpha$ and that $\alpha \neq 0$. 

\medskip

(1) \ Case $d=1$: 
We compute the coefficient $H_{\mathfrak{p}}^{(1)}(\phi)$. 
Since $p(T)= T+ \mu_{0}= T -\alpha$, we get $\mu_{0}= - \alpha$. 
So, we have $\phi_{p(T)}= \phi(T)+ \mu_{0} = p(\alpha)+ A_{1} \tau+ A_{2} \tau^{2}= A_{1} \tau+ A_{2} \tau^{2}$, and hence 
$$
H_{\mathfrak{p}}^{(1)}(\phi)= A_{1}. 
$$

The set $P(1)$ is given by $P(1)= \big\{  (\{ 0 \}, \emptyset)  \big\}$ and $| P(1) |=1$. 

Last, we consider the relation between $H_{\mathfrak{p}}^{(1)}(\phi)$ and $P(1)$. 
The term $A_{1}$ coincides with the element $(\{ 0 \}, \emptyset)$. 
In fact, $A_{1}= L(\emptyset+1) A_{1}^{w(\{0 \})} A_{2}^{w(\emptyset)}$, where $L(\emptyset+1)= 1$ and $w(\emptyset)=0$. 

\medskip

(2) \ Case $d=2$: 
Recall that $\alpha^{q^{2}}= \alpha$. 
We calculate $H_{\mathfrak{p}}^{(2)}(\phi)$. 
Since $p(T)= T^{2}+ \mu_{1} T+ \mu_{0}= (T -\alpha)(T - \alpha^{q})$, we get $\mu_{0}= (-\alpha)(-\alpha^{q})$ and $\mu_{1}= - \alpha- \alpha^{q}$. 
Then, we have 
\begin{align*}
& \phi_{p(T)} = \phi(T)^{2}+ \mu_{1} \phi(T)+ \mu_{0} \\
&= p(\alpha)+ (\alpha^{q}+ \alpha+ \mu_{1})A_{1} \tau + (A_{1}^{q+1}+ \alpha^{q^{2}} A_{2} + \alpha A_{2}+ \mu_{1} A_{2}) \tau^{2} \\
& \qquad + (A_{1} A_{2}^{q}+ A_{1}^{q^{2}} A_{2}) \tau^{3}+ A_{2}^{q^{2}+1} \tau^{4} \\
&=  (A_{1}^{q+1}- (\alpha^{q}- \alpha) A_{2}) \tau^{2} + (A_{1} A_{2}^{q}+ A_{1}^{q^{2}} A_{2}) \tau^{3}+ A_{2}^{q^{2}+1} \tau^{4}, 
\end{align*}
and hence 
$$
H_{\mathfrak{p}}^{(2)}(\phi)= A_{1}^{q+1}- [1] A_{2}= \left( j(\phi)- [1] \right) A_{2},  
$$
where $j(\phi)= A_{1}^{q+1}/A_{2}$. 
Notice that $P_{2}(j(\phi)) \equiv j(\phi)- [1]= H_{\mathfrak{p}}^{(2)}(\phi)/A_{2}$ (mod $\mathfrak{p}$). 

The set $P(2)$ is given by $P(2)= \big\{  (\{ 0,1 \}, \emptyset),  (\emptyset, \{ 0 \})  \big\}$ and $| P(2) |=2$. 

Last, we consider the relation between $H_{\mathfrak{p}}^{(2)}(\phi)$ and $P(2)$. 
The first term $A_{1}^{q+1}$ of $H_{\mathfrak{p}}^{(2)}(\phi)$ coincides with the first element $( \{0, 1 \}, \emptyset)$ of $P(d)$. 
Indeed, $A_{1}^{q+1}= L(\emptyset+1) A_{1}^{w(\{0,1 \})} A_{2}^{w(\emptyset)}$. 
% Recall that $L(\emptyset+1)=L(\emptyset)= 1$ and $w(\emptyset)=0$. 
The second term $-[1] A_{2}$ corresponds to the second element $\{ \emptyset, \{0 \}\}$. 
In fact, $-[1] A_{2}= L(\{1 \}) A_{1}^{w(\emptyset)} A_{2}^{w(\{0\})}$. 
% Recall that $L(\emptyset+1)=L(\emptyset)= 1$ and $w(\emptyset)=0$.   

\medskip

(3) \ Case $d=3$: 
By the same computation as in the above item (2), we obtain 
\begin{align*}
H_{\mathfrak{p}}^{(3)}(\phi) &= A_{1}^{q^{2}+ q+ 1}- [1] A_{1}^{q^{2}} A_{2}- [2] A_{1} A_{2}^{q} \\
& = (j(\phi)^{q}- [1] j(\phi)^{q-1}- [2]) A_{1} A_{2}^{q} \\
& \equiv P_{3}(j(\phi)) A_{1} A_{2}^{q} \pmod {\mathfrak{p}}, 
\end{align*}
% Notice that $P_{3}(j(\phi)) \equiv j(\phi)^{q}- [1]_{\alpha} j(\phi)^{q-1}- [2]_{\alpha}$ (mod $\mathfrak{p}$). 
and $P(3)= \big\{  (\{ 0,1,2 \}, \emptyset),  (\{ 2 \}, \{ 0 \}),  (\{ 0 \}, \{ 1 \})  \big\}$ and $| P(3) |=3$. 

We can check that the $i$th term of $H_{\mathfrak{p}}^{(3)}(\phi)$ coincides exactly with the $i$th element of $P(3)$. 

\medskip

(4) \ Case $d=4$: 
By the same computation as in item (2), we obtain 
\begin{align*}
& H_{\mathfrak{p}}^{(4)}(\phi) = A_{1}^{q^{3}+ q^{2}+q+1} - [1] A_{1}^{q^{3}+ q^{2}} A_{2} - [2] A_{1}^{q^{3}+ 1} A_{2}^{q} \\
& \qquad - [3] A_{1}^{q+ 1} A_{2}^{q^{2}} + [1] [3] A_{2}^{q^{2}+1} \\
& = \Big( j(\phi)^{q^{2+1}} - [1] j(\phi)^{q^{2}} - [2] j(\phi)^{q^{2}-q+1} - [3] j(\phi)+ [1] [3] \Big) A_{2}^{q^{2}+1} \\
& \equiv P_{4}(j(\phi)) A_{2}^{q^{2}+1} \pmod {\mathfrak{p}}, 
\end{align*}
% Notice that $P_{4}(j(\phi)) \equiv j(\phi)^{q^{2+1}} - [1]_{\alpha} j(\phi)^{q^{2}} - [2]_{\alpha} j(\phi)^{q^{2}-q+1} - [3]_{\alpha} j(\phi)+ [1]_{\alpha} [3]_{\alpha}$ (mod $\mathfrak{p}$). 
and 
\begin{align*}
P(4) &= \big\{  (\{ 0,1,2,3 \}, \emptyset),  (\{ 2,3 \}, \{ 0 \}),  (\{ 0,3 \}, \{ 1 \}),  \\
& (\{ 0,1 \}, \{ 2 \}),  (\emptyset, \{ 0,2 \})  \big\}
\end{align*}
and $| P(4) |=5$. 

We can check that the $i$th term of $H_{\mathfrak{p}}^{(4)}(\phi)$ coincides exactly with the $i$th element of $P(4)$. 

\medskip

(5) \ Case $d=5$: 
By the same computation as in item (2), we obtain 
\begin{align*}
& H_{\mathfrak{p}}^{(5)}(\phi) = A_{1}^{q^{4}+ q^{3}+ q^{2}+ q+ 1} - [1] A_{1}^{q^{4}+ q^{3}+ q^{2}} A_{2} - [2] A_{1}^{q^{4}+ q^{3}+ 1} A_{2}^{q} \\
& \qquad - [3] A_{1}^{q^{4}+ q+ 1} A_{2}^{q^{2}} - [4] A_{1}^{q^{2}+ q+ 1} A_{2}^{q^{3}} + [1] [3] A_{1}^{q^{4}} A_{2}^{q^{2}+ 1} \\
& \qquad + [1] [4] A_{1}^{q^{2}} A_{2}^{q^{3}+ 1} + [2] [4] A_{1} A_{2}^{q^{3}+ q} \\
& = \Big( j(\phi)^{q^{3}+q} - [1] j(\phi)^{q^{3}+q-1} - [2] j(\phi)^{q^{3}} - [3] j(\phi)^{q^{3}-q^{2}+q} \\
& \qquad - [4] j(\phi)^{q} + [1] [3] j(\phi)^{q^{3}-q^{2}+q-1} + [1] [4] j(\phi)^{q-1} + [2] [4] \Big) A_{1} A_{2}^{q^{3}+q} \\
& \equiv P_{5}(j(\phi)) A_{1} A_{2}^{q^{3}+q} \pmod {\mathfrak{p}}, 
\end{align*}
% Notice that $P_{5}(j(\phi)) \equiv j(\phi)^{q^{3}+q} - [1]_{\alpha} j(\phi)^{q^{3}+q-1} - [2]_{\alpha} j(\phi)^{q^{3}} - [3]_{\alpha} j(\phi)^{q^{3}-q^{2}+q} - [4]_{\alpha} j(\phi)^{q} + [1]_{\alpha} [3]_{\alpha} j(\phi)^{q^{3}-q^{2}+q-1} + [1]_{\alpha} [4]_{\alpha} j(\phi)^{q-1} + [2]_{\alpha} [4]_{\alpha}$ (mod $\mathfrak{p}$). 
and 
\begin{align*}
P(5) &= \big\{  (\{ 0,1,2,3,4 \}, \emptyset),  (\{ 2,3,4 \}, \{ 0 \}),  (\{ 0,3,4 \}, \{ 1 \}),  (\{ 0,1,4 \}, \{ 2 \}), \\
& (\{ 0,1,2 \}, \{ 3 \}),  (\{ 4 \}, \{ 0,2 \}),  (\{ 2 \}, \{ 0,3 \}),  (\{ 0 \}, \{ 1,3 \})  \big\}
\end{align*}
and $| P(5) |=8$. 

We can check that the $i$th term of $H_{\mathfrak{p}}^{(5)}(\phi)$ coincides exactly with the $i$th element of $P(5)$. 
\end{ex}
%%%%%%%

% \medskip
     
\bigskip

Let $\Z[X_{0}, X_{1}, \ldots, X_{d-1}]$ denote a polynomial ring over $\Z$. 
% Let $\Z[X_{0}, X_{1}, \ldots, X_{d-1}]$ denote the ring of polynomials over $\Z$. 
Let $s_{i} \in \Z[X_{0}, X_{1}, \ldots, X_{d-1}]$ denote the elementary symmetric polynomial of degree $i$ in $-X_{0}$, $-X_{1}$, $\ldots$, $-X_{d-1}$. 
That is, 
\begin{align*}
s_{0} = s_{0}(X_{0}, X_{1}, \ldots, X_{d-1}) & := 1; \\
s_{1} = s_{1}(X_{0}, X_{1}, \ldots, X_{d-1}) & := \sum_{i=0}^{d-1} (-X_{i}); \\
s_{2} = s_{2}(X_{0}, X_{1}, \ldots, X_{d-1}) & := \sum_{0 \leq i<j \leq d-1} (-X_{i}) \cdot (-X_{j}); \\
& \vdots \\
s_{d} = s_{d}(X_{0}, X_{1}, \ldots, X_{d-1}) & := \prod_{i=0}^{d-1} (-X_{i}). 
\end{align*}
For a finite subset $S^{\prime}$ of $\N$ and $n \in \Z$, we define a polynomial $h_{n}^{S^{\prime}}$ in $\Z[ \, X_{i} \ | \ i \in S^{\prime} \, ]$ as follows: 
\begin{align*}
h_{n}^{S^{\prime}} &= h_{n}^{S^{\prime}}( \, X_{i} \ | \ i \in S^{\prime} \, ) \\
& := 
\begin{cases}
\ \displaystyle \sum_{(k_{i}) \in I_{n}(S^{\prime})} \prod_{i \in S^{\prime}} X_{i}^{k_{i}} & \text{if $S^{\prime} \neq \emptyset$}; \\
\ 1 & \text{if $S^{\prime} = \emptyset$ and $n=0$}; \\
\ 0 & \text{if $S^{\prime} = \emptyset$ and $n \neq 0$},  
\end{cases}
\end{align*}
where 
$$
I_{n}(S^{\prime}) := \left\{ (k_{i})_{i \in S^{\prime}} \ \bigg| \ k_{i} \in \N \quad \text{and} \quad \sum_{i \in S^{\prime}} k_{i}= n \right\}. 
$$
Recall that if $n<0$, then $I_{n}(S^{\prime})= \emptyset$, which implies $h_{n}^{S^{\prime}}=0$. 

\medskip

The following plays an very important role in the proof of Main theorem (1). 

%%%%%%%%%%%%%%%%
\begin{Keylemma*}
Assume that $(S_{1}, S_{2}) \in P(d)$, and set $S:= S_{1} \cup S_{2}$ and $S^{\prime}:= S \cup \{ d \}$. 
Then 
\begin{align*} % \label{KEY}
\sum_{i= \lceil d/2 \rceil}^{d} s_{d-i} \cdot h_{i- |S|}^{S^{\prime}} = \prod_{i \in S_{2}} \left( X_{d}- X_{i+1} \right).  
\end{align*}
\end{Keylemma*}
%%%%%%%%%%%%%%

\medskip

%%%%%%%%%%%
\begin{proof}
Let 
$$
L:= \sum_{i= \lceil d/2 \rceil}^{d} s_{d-i} \cdot h_{i- |S|}^{S^{\prime}}, 
$$
and expand this sum into a polynomial. 
Then, each monomial can be uniquely written as 
$$
X_{d}^{k_{d}} \times \prod_{b \in S_{2}} (-X_{b+1})^{\delta_{b}} \times \prod_{c \in S} X_{c}^{k_{c}} \times \prod_{c \in S} (-X_{c})^{\delta_{c}}, 
$$
where $k_{d}$ and $k_{c}$ are nonnegative integers, and $\delta_{b}$ and $\delta_{c}$ are each equal to either $0$ or $1$. 
For simplicity, let 
\begin{align*}
A & := X_{d}^{k_{d}}; \\
B & := \prod_{b \in S_{2}} (-X_{b+1})^{\delta_{b}}; \\
C & := D \cdot E; \\
D & := \prod_{c \in S} X_{c}^{k_{c}}; \\
E & := \prod_{c \in S} (-X_{c})^{\delta_{c}}. 
\end{align*}
Since $\{ S_{1}, S_{2}, S_{1}+1\}$ is a partition, the terms $A$ and $D$ arise from $h_{i- |S|}^{S^{\prime}}$ in $L$ and not from $s_{d-i}$ in $L$. 
Conversely, the other terms $B$ and $E$ arise from $s_{d-i}$ and not from $h_{i- |S|}^{S^{\prime}}$. 

For each term $ABC$, define 
\begin{align*}
\beta &= \beta_{ABC}:= | \{ b \in S_{2} \ | \ \delta_{b} \neq 0 \} |; \\
\gamma &= \gamma_{ABC}:= | \{ c \in S \ | \ \delta_{c} \neq 0 \} |. 
\end{align*}

From here, we compute the sum of terms $ABC$ in $L$ in two ways: first, where $C=1$, and then where $C \neq 1$. 

\medskip

First, we consider the sum of terms $ABC$ with $C=1$ (and then $k_{c}= \delta_{c}= 0$ for any $c \in S$). 
Recall that  $k_{d}+ \sum_{c \in S} k_{c}= i- |S|$, by the definition of $I_{i-|S|}(S^{\prime})$, and so $k_{d}= i- d+ |S_{2}|$. 
Note that $0 \leq \beta \leq |S_{2}|$. 
First, the term $AB$ with $\beta=0$ is exactly the term $X_{d}^{|S_{2}|}$. 
Second, the terms with $\beta=1$ have the form $(- X_{b+1}) X_{d}^{|S_{2}|-1}$ ($b \in S_{2}$).  
Third, the terms with $\beta=2$ have the form $(- X_{b_{1}+1})(- X_{b_{2}+1}) X_{d}^{|S_{2}|-2}$ ($b_{1}, b_{2} \in S_{2}$, $b_{1} \neq b_{2}$). 
In general, terms with $\beta=n$ have the form $(- X_{b_{1}+1})(- X_{b_{2}+1}) \cdots (- X_{b_{n}+1}) X_{d}^{|S_{2}|-n}$ ($b_{1}, b_{2}, \ldots, b_{n} \in S_{2}$, $b_{i} \neq b_{j}$). 
Last, the term with $\beta= |S_{2}|$ is exactly the term $\prod_{b \in S_{2}} (-X_{b+1})$. 
Hence, the sum from $\beta=0$ to $|S_{2}|$ equals 
\begin{align*}
& X_{d}^{|S_{2}|} + \sum_{b \in S_{2}} (- X_{b+1}) X_{d}^{|S_{2}|-1} \\
& \qquad + \sum_{b_{1}, b_{2} \in S_{2}, \atop b_{1}< b_{2}} (- X_{b_{1}+1})(- X_{b_{2}+1}) X_{d}^{|S_{2}|-2} \\
& \qquad + \cdots+ \prod_{b \in S_{2}} (-X_{b+1}) \\
& = \prod_{i \in S_{2}} \left( X_{d}- X_{i+1} \right), 
\end{align*}
and the computation of the first half is complete. 

\medskip 

Next, we consider the sum of the other terms $ABC$ (that is, where $C \neq 1$), and show that the sum is equal to $0$. 
By using the notation $\beta$ and $\gamma$, each term $ABC$ can be uniquely rewritten as follows: 
\begin{align*}
A B C &= (-1)^{\beta+ \gamma} \times X_{d}^{k_{d}} \cdot \prod_{b \in S_{2}} X_{b+1}^{\delta_{b}} \times \prod_{c \in S} X_{c}^{k_{c}+ \delta_{c}}; \\
F &:= X_{d}^{k_{d}} \cdot \prod_{b \in S_{2}} X_{b+1}^{\delta_{b}}. 
\end{align*}
Let 
\begin{align*}
G &:= \prod_{c \in S} X_{c}^{k_{c}^{\prime}}; \\
N &:= N_{G}:= | \{ c \in S \ | \ k_{c}^{\prime} \neq 0 \} |. 
\end{align*}

Now, fix a term $F$, which actually exists and is in $L$, and with this fixed term $F$, fix a term $G$ (a sequence $(k_{c}^{\prime})_{c \in S}$) such that the term $FG$ actually exists in $L$. 
From here, for this fixed sequence $(k_{c}^{\prime})_{c \in S}$, we count the number of terms $ABC$ such that $(k_{c}+ \delta_{c})_{c \in S}= (k_{c}^{\prime})_{c \in S}$ (that is, such that $\prod_{c \in S} X_{c}^{k_{c}+ \delta_{c}}= G$). 
Note that $0 \leq \gamma \leq N$. 
First, the term $ABC$ with $\gamma=0$ is exactly the term such that $(k_{c})_{c \in S}= (k_{c}^{\prime})_{c \in S}$. 
Second, the terms with $\gamma=1$ are terms with a form such that 
\begin{align*}
(\ldots, k_{c}+1, \ldots) &= (\ldots, k_{c}^{\prime}, \ldots); \\
k_{c}^{\prime} &\neq 0  
\end{align*}
($c \in S$), where the other components equal each other. 
Then, the number of such terms is $\binom{N}{1}$. 
Third, terms with $\gamma= 2$ are with a form such that 
\begin{align*}
(\ldots, k_{c_{1}}+1, \ldots, k_{c_{2}}+1, \ldots) &= (\ldots, k_{c_{1}}^{\prime}, \ldots, k_{c_{2}}^{\prime}, \ldots); \\
k_{c_{1}}^{\prime} \neq 0, & \qquad k_{c_{2}}^{\prime} \neq 0
\end{align*}
($c_{1}, c_{2} \in S$, $c_{1} \neq c_{2}$), where the others equal each other. 
Then, the number of such terms is $\binom{N}{2}$.  
In general, the number of terms with $\gamma=n$ is $\binom{N}{n}$. 
So, the sum from $\gamma=0$ to $N$ equals 
\begin{align*}
\sum_{\gamma=0}^{N} (-1)^{\beta+ \gamma} \binom{N}{\gamma} FG &= (-1)^{\beta} \sum_{\gamma=0}^{N} (-1)^{\gamma} \binom{N}{\gamma} FG \\
&= (-1)^{\beta} \left( 1+ (-1) \right)^{N} FG =0. 
\end{align*}
From this, the sum of the terms $ABC$ such that $C \neq 1$ also equals $0$. 
Keylemma follows from this. 
\end{proof}
%%%%%%%%%

\medskip

Here, we explain the structure of the proof of Keylemma. % using Examples \ref{exampless} (5) and \ref{examplepart} (5). 

%%%%%%%%
\begin{ex}
Suppose that $d=5$ and $(S_{1}, S_{2})= (\{ 0 \}, \{ 1,3 \}) \in P(5)$. 
Let $S= S_{1} \cup S_{2}= \{ 0,1,3 \}$ and $S^{\prime}= S \cup \{ 5 \}= \{ 0,1,3,5 \}$. 
Since the sets $I_{n}(S^{\prime})$ are given by
\begin{align*}
I_{0}(S^{\prime}) &= \big\{  (0,0,0,0)  \big\}; \\
I_{1}(S^{\prime}) &= \big\{  (1,0,0,0), \ (0,1,0,0), \ (0,0,1,0), \ (0,0,0,1) \big\}; \\
I_{2}(S^{\prime}) &= \big\{  (1,1,0,0), \ (1,0,1,0), \ (1,0,0,1), \ (0,1,1,0), \\ 
& \qquad (0,1,0,1), \ (0,0,1,1), \ (2,0,0,0), \ (0,2,0,0), \\
& \qquad (0,0,2,0), \ (0,0,0,2) \big\} 
\end{align*}
by definition, the polynomials $h_{n}^{S^{\prime}}= h_{n}^{S^{\prime}}(X_{0}, X_{1}, X_{3}, X_{5})$ can be written as 
\begin{align*}
h_{0}^{S^{\prime}} &= 1;  \\
h_{1}^{S^{\prime}} &= X_{0}+ X_{1}+ X_{3}+ X_{5};  \\
h_{2}^{S^{\prime}} &= X_{0} X_{1}+ X_{0} X_{3}+ X_{0} X_{5}+ X_{1} X_{3} \\
& \qquad + X_{1} X_{5} + X_{3} X_{5} + X_{0}^{2}+ X_{1}^{2}+ X_{3}^{2}+ X_{5}^{2}. 
\end{align*}
% respectively. 
The polynomials $s_{n}= s_{n}(X_{0}, X_{1}, X_{2}, X_{3}, X_{4})$ are defined by 
\begin{align*}
s_{0} &= 1;  \\
s_{1} &= (-X_{0})+ (-X_{1}) + (-X_{2}) + (-X_{3}) + (-X_{4});  \\
s_{2} &= (-X_{0}) (-X_{1}) + (-X_{0}) (-X_{2}) + (-X_{0}) (-X_{3}) + (-X_{0}) (-X_{4}) \\
& \qquad + (-X_{1}) (-X_{2}) + (-X_{1}) (-X_{3}) + (-X_{1}) (-X_{4}) \\
& \qquad + (-X_{2}) (-X_{3}) + (-X_{2}) (-X_{4}) + (-X_{3}) (-X_{4}). 
\end{align*}
Therefore, we obtain 
\begin{align*}
& \sum_{i= \lceil 5/2 \rceil}^{5} s_{5-i} \cdot h_{i- 3}^{S^{\prime}} = s_{2} \cdot h_{0}^{S^{\prime}}+ s_{1} \cdot h_{1}^{S^{\prime}}+ s_{0} \cdot h_{2}^{S^{\prime}} \\
&= \Big(  (-X_{0})(-X_{1}) + (-X_{0})(-X_{2}) + (-X_{0})(-X_{3}) + (-X_{0})(-X_{4}) \\
& \qquad + (-X_{1})(-X_{2}) + (-X_{1})(-X_{3}) + (-X_{1})(-X_{4}) \\
& \qquad + (-X_{2})(-X_{3}) + (-X_{2})(-X_{4}) + (-X_{3})(-X_{4})  \Big) \cdot 1 \\
& + \Big(  (-X_{0})+ (-X_{1}) + (-X_{2}) + (-X_{3}) + (-X_{4})  \Big) \\
& \qquad \times \Big(  X_{0} + X_{1} + X_{3} + X_{5}  \Big) \\
& + 1 \cdot  \Big(  X_{0} X_{1} + X_{0} X_{3} + X_{0} X_{5} + X_{1} X_{3} + X_{1} X_{5} \\
& \qquad + X_{3} X_{5} + X_{0}^{2} + X_{1}^{2} + X_{3}^{2} + X_{5}^{2}  \Big), 
\end{align*}
and hence 
\begin{align*} 
& \sum_{i= \lceil 5/2 \rceil}^{5} s_{5-i} \cdot h_{i- 3}^{S^{\prime}} = X_{5}^{2} + \Big(  (-X_{2}) + (-X_{4})  \Big) X_{5} + (-X_{2}) (-X_{4}) \\
& \qquad + \Big(  1 + 2 (-1) + (-1)^{2}  \Big)  \Big(  X_{0} X_{1} + X_{0} X_{3} + X_{1} X_{3}  \Big) \\
& \qquad + (-1)  \Big(  1 + (-1)  \Big)  \Big(  X_{0} X_{2} + X_{0} X_{4} + X_{1} X_{2} + X_{1} X_{4} + X_{3} X_{4}  \Big) \\
& \qquad + \Big(  1 + (-1)  \Big)  \Big(  X_{0} X_{5}+ X_{1} X_{5}+ X_{3} X_{5} \Big) \\
& \qquad + \Big(  1 + (-1)  \Big)  \Big(  X_{0}^{2}+ X_{1}^{2}+ X_{3}^{2}  \Big) \\
&= (X_{5}- X_{2})(X_{5}- X_{4}).     
\end{align*}

We now consider the structure in detail. 
The terms $X_{5}^{2}$, $(-X_{2}) X_{5}$, $(-X_{4}) X_{5}$, and $(-X_{2}) (-X_{4})$ are the terms $ABC$ such that $C=1$ in the proof of Keylemma. 
In contrast, the other terms are the terms $ABC$ such that $C \neq 1$. 
Now, we fix the terms $F=1$ and $G= X_{0} X_{1}$ (that is, the sequence $(k_{0}^{\prime}, k_{1}^{\prime}, k_{3}^{\prime})= (1,1,0)$). 
Note that $\beta=0$ and $N=2$. 
We count the number of terms such that $(k_{0}+ \delta_{0}, k_{1}+ \delta_{1}, k_{3}+ \delta_{3})= (1,1,0)$. 
First, the term with $\gamma=0$ is exactly the term such that $(1+0,1+0,0)= (1,1,0)$. 
Second, the terms with $\gamma= 1$ are the terms such that either $(0+1,1+0,0)= (1,1,0)$ or $(1+0, 0+1,0)= (1,1,0)$. 
Last, the term with $\gamma= 2$ is exactly the term such that $(0+1,0+1,0)= (1,1,0)$. 
Hence, the sum equals 
$$
X_{0} X_{1}+ 2(-1) X_{0} X_{1}+ (-1)^{2} X_{0} X_{1}= 0. 
$$
Similarly, the terms $X_{0} X_{3}$ and $X_{1} X_{3}$ can be computed. 

Next, we fix the terms $F=X_{2}$ and $G= X_{0}$ (that is, the sequence $(k_{0}^{\prime}, k_{1}^{\prime}, k_{3}^{\prime})= (1,0,0)$). 
Note that $\beta=1$ and $N=1$. 
We count the number of terms such that $(k_{0}+ \delta_{0}, k_{1}+ \delta_{1}, k_{3}+ \delta_{3})= (1,0,0)$. 
Then, the term with $\gamma=0$ (resp. $\gamma=1$) is exactly the term such that $(1+0,0,0)= (1,0,0)$ (resp. $(0+1,0,0)= (1,0,0)$). 
Hence, the sum equals 
$$
(-1) X_{0} X_{2}+ (-1)^{1+1} X_{0} X_{2}=0. 
$$
Similarly, the terms $X_{0} X_{4}$, $X_{1} X_{2}$, $X_{1} X_{4}$, $X_{3} X_{4}$, $X_{0} X_{5}$, $X_{1} X_{5}$, and $X_{3} X_{5}$ can be calculated. 

Finally, we fix the terms $F=1$ and $G= X_{0}^{2}$ (that is, the sequence $(k_{0}^{\prime}, k_{1}^{\prime}, k_{3}^{\prime})= (2,0,0)$). 
Note that $\beta=0$ and $N=1$. 
We count the number of terms such that $(k_{0}+ \delta_{0}, k_{1}+ \delta_{1}, k_{3}+ \delta_{3})= (2,0,0)$. 
First, the term with $\gamma=0$ (resp. $\gamma=1$) is the term such that $(2+0,0,0)= (2,0,0)$ (resp. $(1+1,0,0)= (2,0,0)$). 
Hence, the sum equals 
$$
X_{0}^{2}+ (-1) X_{0}^{2}=0. 
$$
Similarly, the terms $X_{1}^{2}$ and $X_{3}^{2}$ can be computed, which finishes the illustration of Keylemma. 
\end{ex}
%%%%%%%

\bigskip

% \medskip

Recall that $\alpha, \alpha^{q}, \ldots, \alpha^{q^{d-1}}$ are the roots of $p(T)= \sum_{i=0}^{d} \mu_{i} T^{i}$, and that $\alpha^{q^{d}}= \alpha$. 
% Recall that $\alpha \neq 0$. 
% In fact, $p(\alpha^{q^{i}})= p(\alpha)^{q^{i}}=0$. 
% In fact, $p(\alpha^{q^{i}})= p(\alpha)^{q^{i}}=0$. 
% Then, $A/\mathfrak{p}= \F_{q^{d}}= \F_{q}(\alpha)$, and 
By the relation between the roots of $p(T)$ and the coefficients of $p(T)$, we get 
\begin{align*}
\mu_{d} &:= 1; \\
\mu_{d-1} &= \sum_{i=0}^{d-1} (- \alpha^{q^{i}}); \\
\mu_{d-2} &= \sum_{0 \leq i<j \leq d-1} (- \alpha^{q^{i}})(- \alpha^{q^{j}}); \\
& \vdots \\
\mu_{0} &= \prod_{i=0}^{d-1} (- \alpha^{q^{i}}). 
\end{align*}
Then, it follows from the definitions of $s_{d-i}$ and $\mu_{i}$ that $s_{d-i}(\alpha, \alpha^{q}, \ldots, \alpha^{q^{d-1}}) = \mu_{i}$. % and $\alpha^{q^{d}}= \alpha$. 

\bigskip

% \medskip

As applications of Keylemma, we have the following corollary.    

%%%%%%%%%
\begin{cor}\label{corkeylem}
% Let $d \geq 1$ be an integer. 
Assume that $(S_{1}, S_{2}) \in P(d)$, and let $S:= S_{1} \cup S_{2}$ and $S^{\prime}:= S \cup \{ d \}$. 

\medskip

(a) \ 
Further, suppose that $L$ is an extension of $K$. 
Then 
\begin{align*}
\sum_{i= \lceil d/2 \rceil}^{d} s_{d-i}(T, T^{q}, \ldots, T^{q^{d-1}}) & \cdot h_{i- |S|}^{S^{\prime}} ( \, T^{q^{i}} \ | \ i \in S \, ) \\
& = \prod_{i \in S_{2}} \left( [d]- [i+1] \right). 
\end{align*}

\medskip

(b) \ 
Instead, assume that $L$ is an extension of $\F_{\mathfrak{p}}$.
Then  
% and moreover, if we substitute $T= \alpha$ in the above claim, then we obtain 
\begin{align*}
\sum_{i= \lceil d/2 \rceil}^{d} \mu_{i} \cdot h_{i- |S|}^{S^{\prime}} ( \, \alpha^{q^{i}} \ | \ i \in S \, ) & = (-1)^{|S_{2}|} \prod_{i \in S_{2}} [i+1] \\
& = L(S_{2}+1). 
\end{align*}
\end{cor}
%%%%%%%%

\medskip

%%%%%%%%%%%
\begin{proof}
(a) \ 
Recall that $[i]= T^{q^{i}}- T$. 
If we substitute $X_{i} = T^{q^{i}}$ in Keylemma, then we get $X_{d}- X_{i+1}= T^{q^{d}}- T^{q^{i+1}}= [d]- [i+1]$. 
This proves the first claim. 
% Hence, the first claim is proven. 

\medskip

(b) \ 
Recall that $[d]= \alpha^{q^{d}}- \alpha= 0$. 
The second claim follows from item (a). 
\end{proof}
%%%%%%%%%

\medskip

Let $\phi$ be a rank-$2$ Drinfeld module over $L$, and let $\Lambda_{\phi}$ be the lattice corresponding to $\phi$ (see Theorem 2.4 of Chapter I in \cite{Gek86}, Theorem 4.6.9 in \cite{Gos96}, Theorem 2.4.2 in \cite{Tha04}). 
Let $e_{\phi}(z)$ denote the lattice exponential function of $\Lambda_{\phi}$, defined by $e_{\phi}(z):= z \prod_{0 \neq \lambda \in \Lambda_{\phi}} (1- z/\lambda)$ (see Definition 2.1 of Chapter I in \cite{Gek86}, Definition 4.2.3 in \cite{Gos96}, Subsection 2.4 in \cite{Tha04}).  % ($z \in \C_{\infty}$) 
Then, it is known that $e_{\phi}(z)$ has the composition inverse function $\log_{\phi}(z)$ such that % a series expansion 
$$
\log_{\phi}(z)= \sum_{j \geq 0} \beta_{j} z^{q^{j}}
$$
(see Section 2 of Chapter II in \cite{Gek86}, Subsection 4.6 in \cite{Gos96}, Subsection 2.4 in \cite{Tha04}). 
Note that $\beta_{0}=1$ from the definition of $e_{\phi}(z)$. 

\medskip

We now collect results of El-Guindy and Papanikolas to use in the proof of Main theorem (1). 

%%%%%%%%%%%
\begin{fact*}[Theorem 8.1, Corollary 8.2 and Theorem 3.3 in \cite{EP13}]
% Let $\mathfrak{p}= (p(T))$ be a nonzero prime ideal of $A$, where $p(T) = \sum_{i=0}^{d} \mu_{i} T^{i}$ is a monic irreducible polynomial of degree $d$. 
% Assume that $L$ is an extension of either $K$ or $\F_{\mathfrak{p}}$. 
Let $\phi$ be any rank-$2$ Drinfeld module over a field $L$ defined by $\phi_{T}= \iota(T) + A_{1} \tau + A_{2} \tau^{2}$. 

\medskip

(1) \ 
% Assume that $L=K$. 
Let $m \in \N$ and $n \in \Z$. 
Now, we define coefficients $c(n; m):= c(n; m; \phi)$ as follows. 
For $0 \leq n \leq 2m$, let 
$$
\phi_{T^{m}}= \sum_{n=0}^{2m} c(n; m) \tau^{n}. 
$$
For the other cases ($n<0$ or $n> 2m$), set $c(n; m)=0$. 
Then, for any $m,n \geq 0$, we have that 
$$
c(n; m)= \sum_{(S_{1}, S_{2}) \in P(n)} A_{1}^{w(S_{1})} A_{2}^{w(S_{2})} \cdot h_{m- |S|}^{S^{\prime}}, 
$$
where $S := S_{1} \cup S_{2}$ and $S^{\prime}:= S \cup \{ n \}$. 

\medskip

(2) \ 
Next, suppose that $L$ is an extension of $\F_{\mathfrak{p}}$. 
Then 
$$
H_{\mathfrak{p}}^{(d)}(\phi) = \sum_{i= \lceil d/2 \rceil}^{d} \mu_{i} \cdot c(d; i). 
$$
Moreover, $\phi$ is supersingular at $\mathfrak{p}$ if and only if 
$$
\sum_{i= \lceil d/2 \rceil}^{d} \mu_{i} \cdot c(d; i) = 0. 
$$

\medskip

(3) \ 
The coefficients $\beta_{j}$ of the function $\log_{\phi}(z)$ are given by 
$$
\beta_{j}= \sum_{(S_{1}, S_{2}) \in P(j)} \frac{A_{1}^{w(S_{1})} A_{2}^{w(S_{2})}}{L(S_{1}+1) L(S_{2}+2)}. 
$$
\end{fact*}
%%%%%%%%%

\medskip

% Let $d \geq 1$ be a positive integer. 
Let $L$ be an extension of $K$, and let $\phi$ be a rank-$2$ Drinfeld module over $L$, defined by $\phi_{T}= T + A_{1} \tau + A_{2} \tau^{2}$. 
Recall that $\delta^{q}= - \delta$ (see Section \ref{sectintro}). 
Set 
\begin{align*}
\mathfrak{a}(d) := 
\begin{cases}
\ T \cdot \sum_{j=0}^{d} \beta_{j} \delta^{q^{j}} & \text{if $d \geq 1$}; \\
\ T \delta & \text{if $d= 0$}; \\  
\ 0 & \text{if $d<0$} \\  
\end{cases}
\end{align*}
(see Theorem 6.3 in \cite{EP13}). 
It follows from Fact (3) that 
\begin{align}\label{fact3}
\mathfrak{a}(d) &= T \cdot \sum_{j=0}^{d} \beta_{j} \delta^{q^{j}}= T \delta \cdot \sum_{j=0}^{d} (-1)^{j} \beta_{j} \nonumber \\
&= T \delta \cdot \sum_{j=0}^{d} (-1)^{j} \sum_{(S_{1}, S_{2}) \in P(j)} \frac{A_{1}^{w(S_{1})} A_{2}^{w(S_{2})}}{L(S_{1}+1) L(S_{2}+2)}, 
\end{align}
which is used in proving the proposition below. 

For this, we write 
$$
L_{d}:= (-1)^{d} [d][d-1] \cdots [2][1] \qquad \text{and} \qquad L_{0}:= 1. 
$$
We define 
\begin{align*}
\mathfrak{b}(d) := 
\begin{cases}
\ \displaystyle \frac{L_{d} \cdot \mathfrak{a}(d)}{\delta T^{1+ q+ \cdots+ q^{d}}} & \text{if $d \geq 1$}; \\
\ 1 & \text{if $d= 0$}; \\
\ 0 & \text{if $d<0$} \\  
\end{cases}
\end{align*}
(cf. Section 2 in \cite{ElG14}). 

Assume that $A_{1}= T+ \lambda$ and $A_{2}= \lambda$. 
Then, the following recursive relations are valid: 
\begin{align}\label{recursionb}
-[d] \mathfrak{a}(d) &= - \mathfrak{a}(d-1) (T^{q^{d}}+ \lambda^{q^{d-1}})+ \lambda^{q^{d-1}} \mathfrak{a}(d-2); \nonumber \\
\mathfrak{b}(d) &= -(1+ D^{q^{d-1}})  \mathfrak{b}(d-1)+ D^{q^{d-1}} (T^{1- q^{d-1}}-1) \mathfrak{b}(d-2), 
\end{align}
where $D:= {\lambda}/{T^{q}}$. 
The validity of these recursions can be proven in the same way as for the recursions (12) and (14), respectively, in \cite{ElG14}. 
% We shall remark that these recursion formulae are different from Recursion formula (2.1) in \cite{Cor99}, that is, these recursion formulae can not be induced from Recursion formula (2.1). 

\bigskip

% \medskip

The following is Main theorem (1). 

%%%%%%%%%%
\begin{prop}\label{mainthm1}
Let $\mathfrak{p}= (p(T))$ denote a nonzero prime ideal of degree $d$ such that $p(T) \neq T$, and let $\phi$ be any rank-$2$ Drinfeld module over a field $L$ defined by $\phi_{T}= \iota(T)+ A_{1} \tau+ A_{2} \tau^{2}$. 
% Let $\mathfrak{p}= (p(T))$ be a nonzero prime ideal of $A$, where $p(T)= \sum_{i=0}^{d} \mu_{i} T^{i}$ is a monic irreducible polynomial of degree $d$. 
% Assume that $L$ is an extension of either $K$ or $\F_{\mathfrak{p}}$. 
% Assume that $p(T) \neq T$. 

\medskip

(a) \ 
Suppose that $L$ is an extension of $\F_{\mathfrak{p}}$. 
Then, the coefficient $H_{\mathfrak{p}}^{(d)}(\phi)$ is given by 
\begin{align*}
H_{\mathfrak{p}}^{(d)}(\phi) &= \sum_{(S_{1}, S_{2}) \in P(d)} L(S_{2}+1) A_{1}^{w(S_{1})} A_{2}^{w(S_{2})}. % \pmod{\mathfrak{p}}. \\
\end{align*}

\medskip

(b) \ 
Assume that $L$ is an extension of $K$. 
Then 
\begin{align*}
H_{\mathfrak{p}}^{(d)}(\phi) & \equiv \frac{(-1)^{d} L_{d} \cdot \mathfrak{a}(d)}{\delta T^{q^{d}}} \pmod{\mathfrak{p}} \\
& = (-1)^{d} T^{1+q+ \cdots + q^{d-1}} \mathfrak{b}(d). 
\end{align*}
\end{prop}
%%%%%%%%%

\medskip

%%%%%%%%%%%
\begin{proof}
(a) \ 
It follows from Fact (1) and the first claim of Fact (2) that 
\begin{align*}
H_{\mathfrak{p}}^{(d)}(\phi) &= \sum_{i= \lceil d/2 \rceil}^{d} \mu_{i} \cdot c(d; i) \\
&= \sum_{(S_{1}, S_{2}) \in P(d)} \left( \sum_{i= \lceil d/2 \rceil}^{d} \mu_{i} \cdot h_{i- |S|}^{S^{\prime}} \right) A_{1}^{w(S_{1})} A_{2}^{w(S_{2})}. 
\end{align*}
By using Corollary \ref{corkeylem} (b), we have 
\begin{align*}
H_{\mathfrak{p}}^{(d)}(\phi) &= \sum_{(S_{1}, S_{2}) \in P(d)} \left( \sum_{i= \lceil d/2 \rceil}^{d} \mu_{i} \cdot h_{i- |S|}^{S^{\prime}} \right) A_{1}^{w(S_{1})} A_{2}^{w(S_{2})} \\
& = \sum_{(S_{1}, S_{2}) \in P(d)} L(S_{2}+1) A_{1}^{w(S_{1})} A_{2}^{w(S_{2})}. % \pmod{\mathfrak{p}}. \\
\end{align*}
This completes the proof of this case. 

\medskip

(b) \ 
Let $0 \leq j \leq d$ be an integer. 
If $(S_{1}, S_{2}) \in P(j)$, then $\{ S_{1}, S_{2}, S_{2}+1 \}$ is a partition of $\N_{<j}$, and thus $\{ S_{1}+1, S_{2}+1, S_{2}+2 \}$ is a partition of $\{ 1,2, \ldots, j \}$. 
Hence, for $0 \leq j < d$, we obtain 
$$
\frac{L_{d}}{L(S_{1}+1) L(S_{2}+2)}= (-1)^{d-j} [d] [d-1] \cdots [j+1] L(S_{2}+1) \equiv 0 \pmod{\mathfrak{p}}. 
$$
Along another line, for $(S_{1}, S_{2}) \in P(d)$, we get $L_{d}= L(S_{1}+1) L(S_{2}+1) L(S_{2}+2)$. 
Hence, 
\begin{align*}
& \frac{(-1)^{d} L_{d} \cdot \mathfrak{a}(d)}{\delta T^{q^{d}}} \\
&= \frac{(-1)^{d} L_{d}}{\delta T^{q^{d}}} \cdot T \delta \cdot \sum_{j=0}^{d} (-1)^{j} \sum_{(S_{1}, S_{2}) \in P(j)} \frac{A_{1}^{w(S_{1})} A_{2}^{w(S_{2})}}{L(S_{1}+1) L(S_{2}+2)} \\
& \equiv \sum_{(S_{1}, S_{2}) \in P(d)} L(S_{2}+1) A_{1}^{w(S_{1})} A_{2}^{w(S_{2})} = H_{\mathfrak{p}}^{(d)}(\phi) \pmod{\mathfrak{p}}
\end{align*}
by the equality (\ref{fact3}) and item (a). 
The second claim follows. 
\end{proof}
%%%%%%%%%

\medskip

We next consider relations between the known polynomials and our polynomials. 

%%%%%%%%%%%
\begin{rem*}
Assume that $L$ is an extension of $K$, and that $A_{1}= T+ \lambda$ and $A_{2}= \lambda$. 
% Set $D:= \lambda/T^{q}$. 

\medskip

(1) \ 
We explore the relation between the polynomial $\mathfrak{b}_{d}$ in Section 2 of \cite{ElG14} and our polynomial $\mathfrak{b}(d)$. 
First, we easily see that $\mathfrak{b}_{d}= (-1)^{d} \mathfrak{b}(d)$. 
Next, from Theorem 2.2 in \cite{ElG14}, we have 
$$
\mathfrak{b}(d)= (-1)^{d} \sum_{S \subseteq \N_{<d}} \frac{(\lambda/T^{q})^{w(S)}}{m(S)},  
$$
where $M(S):= S \setminus (S+1)$ and $m(S):= \prod_{i \in M(S)} T^{q^{i}-1}$. 

\medskip

(2) \ 
We investigate the relation between the polynomial $p_{d}(s)$ in Definition 12 of \cite{BB} and our polynomial $H_{\mathfrak{p}}^{(d)}(\lambda)$. 
We see that 
\begin{align*}
H_{\mathfrak{p}}^{(d)}(\lambda) \equiv (-1)^{d} T^{1+q+ \cdots + q^{d-1}} p_{d} \left( -\lambda/T^{q} \right) \pmod {\mathfrak{p}} 
\end{align*}
by Proposition \ref{mainthm1} (b), item (1) of this remark, and Theorem 3.2 (i) in \cite{ElG14}. 
% $$ \mathfrak{b}(n)= ((-D)^{q^{n-1}}-1)  \mathfrak{b}(n-1)+ (-D)^{q^{n-1}} (1- T^{1- q^{n-1}}) \mathfrak{b}(n-2) $$
\end{rem*}
%%%%%%%%%

\medskip

Last, we study the relation between one period of a lattice $\Lambda_{\phi}$ and our coefficient $H_{\mathfrak{p}}^{(d)}(\phi)$. 

%%%%%%%%%%
\begin{rem*}
(1) \ 
Assume that $L$ is an extension of $K$. 
It is known that the series 
$$
\mathfrak{f}(z):= \sum_{n \geq 0} \mathfrak{a}(n) z^{q^{n}}
$$
forms a period of $\Lambda_{\phi}$ (see Theorem 6.3 in \cite{EP13}). 
Let $\mathfrak{f}^{\text{tr}}(z)$ denote the truncated function defined by 
$$
\mathfrak{f}^{\text{tr}}(z):= \sum_{n= 0}^{d} \mathfrak{a}(n) z^{q^{n}}. 
$$
Recall that $\delta^{q}= - \delta$ and $\delta^{q^{n}}= (-1)^{n} \delta$. 
From Proposition \ref{mainthm1} (b), we get 
\begin{align*}
L_{d} \cdot \mathfrak{f}^{\text{tr}} \left( \frac{1}{\delta T} \right) &= \sum_{n= 0}^{d} L_{d} \cdot \mathfrak{a}(n) \left( \frac{1}{\delta T} \right)^{q^{n}} \\
&= \sum_{n=0}^{d} \frac{L_{d}}{\delta^{q^{n}} T^{q^{n}}} \frac{\delta T^{1+ q+ \cdots+ q^{n}}}{L_{n}} \mathfrak{b}(n) \\
&= (-1)^{d} \sum_{n=0}^{d} [d][d-1] \cdots [n+1] \frac{T^{1+ q+ \cdots+ q^{n}}}{T^{q^{n}}} \mathfrak{b}(n),   
\end{align*}
where $L_{d}:= (-1)^{d} [d][d-1] \cdots [2][1]$. 
Note that $[d] \equiv 0$ $(\text{mod} \ \mathfrak{p})$. 
By Proposition \ref{mainthm1} (b), we obtain 
\begin{align*}
L_{d} \cdot \mathfrak{f}^{\text{tr}} \left( \frac{1}{\delta T} \right) & \equiv (-1)^{d} T^{1+ q+ \cdots+ q^{d-1}} \mathfrak{b}(d) \\
& \equiv H_{\mathfrak{p}}^{(d)}(\phi) \pmod {\mathfrak{p}}. 
\end{align*}
Hence, we can regard the constant $1/L_{d}$ as the constant $\pi^{\text{tr}}$ introduced in Section \ref{sectintro} by the following reasoning. 

Let 
\begin{align*}
\xi_{\ast} &:= \prod_{n \geq 0} \left( 1- \frac{[n]}{[n+1]} \right); \\
\xi_{d} &:= \frac{[1]^{(q^{d}-1)/(q-1)}}{(-1)^{d} L_{d}}= \prod_{n=0}^{d-1} \left( 1- \frac{[n]}{[n+1]} \right)
% \xi &:= L_{1}^{1/(q-1)} \xi_{\ast}, 
\end{align*}
(see Theorem 1.4 of Chapter IV in \cite{Gek86}, Subsection 3.2 in \cite{Gos96}, Subsection 2.5 in \cite{Tha04}). 
% where $L_{1}^{1/(q-1)}$ is any $(q-1)$-st root of $L_{1}$ (see Definition 3.2.7 in \cite{Gos96}). 
It is known that the constant $\xi_{\ast}$ is a function-field analogue of the circular constant $\pi$ (see Remark 1.5 of Chapter IV in \cite{Gek86}). 
Hence, we can treat the constant $1/L_{d}= (-1)^{d} \xi_{d}/[1]^{(q^{d}-1)/(q-1)}$ as the constant $\pi^{\text{tr}}$. % (???)

\medskip

(2) \ 
The Deuring polynomial $H_{p}(\lambda)$ is related to a hypergeometric function (see Section \ref{sectintro}). 
Thakur defined several hypergeometric functions for function fields (see \cite{Tha95, Tha00}, Subsection 6.5 in \cite{Tha04}). 
However, it seems that the functions given by Thakur are not related to the coefficient $H_{\mathfrak{p}}^{(d)}(\phi)$. 
\end{rem*}
%%%%%%%%%

%%%%%%%%%%%%%%%%%%%%%%%%%%%%%%%%%%%%%%%
\section{Properties for a supersingular polynomial}\label{applications}
%%%%%%%%%%%%%%%%%%%%%%%%%%%%%%%%%%%%%%%

In this section, we prove Main theorem (2) (Proposition \ref{mainthm23} (b), (c), (d)), which was introduced in Section \ref{sectintro}. % , which is a function field analogue of several properties for the Deuring's polynomial. % in the elliptic case. 
In the course of the proof, we often use Main theorem (1). % , which determines the explicit formula. % for the supersingular polynomial. 

\bigskip

% \medskip

We recall the following notation. 
In what follows, $\mathfrak{p}= (p(T))$ denotes a nonzero prime ideal of degree $d$ in $A= \F_{q}[T]$ such that $p(T) \neq T$, and $\alpha$ is a root of $p(T)$. 
Note that $\alpha \neq 0$. 
% Set $d= \deg_{T} p(T)$. 
% Assume that $L$ is an extension of $\F_{\mathfrak{p}}$. 
Throughout this section, we assume that $L$ is an extension of $\F_{\mathfrak{p}}$, and that $\phi$ is a rank-$2$ Drinfeld module over $L$ defined by 
$$
\phi_{T}= \alpha+ (\alpha+ \lambda) \tau+ \lambda \tau^{2}. 
$$
Set 
$$
\phi_{\mathfrak{p}}:= \phi_{p(T)}= \sum_{i=d}^{2d} g_{i} \tau^{i} \qquad \text{and} \qquad H_{\mathfrak{p}}^{(d)}(\lambda):= g_{d}. 
$$

\medskip

The following proposition is Main theorem (2). 

%%%%%%%%%%
\begin{prop}\label{mainthm23} % [cf. Section 5 in \cite{Gek83}]
% The followings hold:  

% \medskip

(a) \ 
The polynomial $H_{\mathfrak{p}}^{(d)}(\lambda)$ is a divisor of $g_{i}$ for all $d \leq i < 2d$. 
Moreover, if $\phi$ is supersingular at $\mathfrak{p}$, then $\phi_{\mathfrak{p}}= \lambda^{(q^{2d}-1)/(q^{2}-1)} \tau^{2d}$. 

\medskip

(b) \ 
$H_{\mathfrak{p}}^{(d)}(0) \neq 0$ and $\deg_{\lambda} H_{\mathfrak{p}}^{(d)}(\lambda)= (q^{d}-1)/(q-1)$.  

\medskip

(c) \ 
$H_{\mathfrak{p}}^{(d)}(\lambda)$ and $H_{\mathfrak{p}}^{(d)}(-\alpha^{q} s(s+1)^{q-1})$ are each separable. % BB, simple

\medskip

(d) \ 
If $H_{\mathfrak{p}}^{(d)}(\lambda^{q+1})=0$, then $\lambda \in \F_{\mathfrak{p}}^{(2)}$. 
In particular, all the roots of $H_{\mathfrak{p}}^{(d)}(\lambda)$ are in $\F_{\mathfrak{p}}^{(2)}$. % $(q+1)$-st, BB?

\end{prop}
%%%%%%%%

\medskip

%%%%%%%%%%%
\begin{proof}
(a) \ 
Note that $A= \F_{q}[T]$ is commutative. 
Recall that $\phi_{T} = \phi(T)= \alpha+ A_{1} \tau+ A_{2} \tau^{2}$ and $\phi_{p(T)} = \phi(p(T))=  \sum_{i=0}^{2d} g_{i} \tau^{i}$. 
Since $\phi$ is a homomorphism, we have 
$$
\phi_{T} \phi_{p(T)}= \phi_{T \cdot p(T)}= \phi_{p(T) \cdot T}= \phi_{p(T)} \phi_{T}. 
$$
Now, we compute the left-hand side and the right-hand side:  
\begin{align*}
\phi_{T} \phi_{p(T)} &= (\alpha+ A_{1} \tau+ A_{2} \tau^{2}) \cdot \sum_{i=0}^{2d} g_{i} \tau^{i} \\
&= \cdots+ (g_{i-2}^{q^{2}} A_{2}+ g_{i-1}^{q} A_{1}+ g_{i} \alpha) \tau^{i} + \cdots, \quad \text{and} \\
\phi_{p(T)} \phi_{T} &= \sum_{i=0}^{2d} g_{i} \tau^{i} \cdot (\alpha+ A_{1} \tau+ A_{2} \tau^{2}) \\
&= \cdots+ (g_{i-2} A_{2}^{q^{i-2}} + g_{i-1} A_{1}^{q^{i-1}} + g_{i} \alpha^{q^{i}}) \tau^{i} + \cdots. 
\end{align*}
Then, the coefficients of $\tau^{i}$ are the same between sides, as shown by 
\begin{align*}
g_{i-2}^{q^{2}} A_{2}+ g_{i-1}^{q} A_{1}+ g_{i} \alpha= g_{i-2} A_{2}^{q^{i-2}} + g_{i-1} A_{1}^{q^{i-1}} + g_{i} \alpha^{q^{i}}, 
\end{align*}
and so 
\begin{align}\label{recurrencegi}
(\alpha^{q^{i}}- \alpha) g_{i}= g_{i-2}^{q^{2}} A_{2}- g_{i-2} A_{2}^{q^{i-2}}+ g_{i-1}^{q} A_{1}- g_{i-1} A_{1}^{q^{i-1}}. 
\end{align}
Recall that $g_{d-1}= 0$ and $g_{d}= H_{\mathfrak{p}}^{(d)}(\lambda)$. 
Note that $\alpha^{q^{i}}- \alpha \neq 0$ for all $d< i< 2d$. 
First, suppose that $i=d+1$ for the recursion (\ref{recurrencegi}). 
Since $g_{d-1}=0$, we get $g_{d} \mid g_{d+1}$. 
Second, assume that $i=d+2$. 
Since $g_{d} \mid g_{d+1}$, we get $g_{d} \mid g_{d+2}$. 
We obtain $g_{d} \mid g_{i}$ for all $d \leq i < 2d$ by induction on $i$. 

\medskip

Recall that $p(T)= \sum_{i=0}^{d} \mu_{i} T^{i}$ and $A_{2}= \lambda$.  
Note that $g_{0}= \cdots= g_{d-1}= 0$. 
Notice that 
\begin{align*}
\phi_{T^{2}} &=  \phi_{T} \phi_{T} \\
&= (\alpha+ A_{1} \tau+ A_{2} \tau^{2})  (\alpha+ A_{1} \tau+ A_{2} \tau^{2})= \cdots+ A_{2}^{q^{2}+1} \tau^{4}, \quad \text{and} \\
\phi_{T^{3}} &= \phi_{T} \phi_{T^{2}} \\
&= (\alpha+ A_{1} \tau+ A_{2} \tau^{2}) (\cdots+ A_{2}^{q^{2}+1} \tau^{4})= \cdots+ A_{2}^{q^{4}+ q^{2}+ 1} \tau^{6}. 
\end{align*}
Hence, we have $\phi_{T^{i}}= \cdots+ A_{2}^{q^{2(i-1)}+ \cdots+ q^{2}+ 1} \tau^{2i}$ for all $1 \leq i \leq d$, again by induction on $i$. 
From this, we obtain $\phi_{p(T)} = \sum_{i=0}^{d} \mu_{i} \phi_{{T}^{i}}= \cdots+ A_{2}^{q^{2(d-1)}+ \cdots+ q^{2}+ 1} \tau^{2d}$, namely, $g_{2d}= \lambda^{q^{2(d-1)}+ \cdots+ q^{2}+ 1}$. 

If $\phi$ is supersingular at $\mathfrak{p}$, then we get $g_{d}= H_{\mathfrak{p}}^{(d)}(\lambda)=0$ (see Satz (5.3) in \cite{Gek83}). 
It follows from the first claim that $g_{d+1}= \cdots= g_{2d-1}=0$. 
Thus, we have $\phi_{\mathfrak{p}}= \phi_{p(T)}= \lambda^{q^{2(d-1)}+ \cdots+ q^{2}+ 1} \tau^{2d}$. 

\medskip

(b) \ 
By Main theorem (1), we have 
$$
H_{\mathfrak{p}}^{(d)}(\lambda)= \sum_{(S_{1}, S_{2}) \in P(d)} L(S_{2}+1) (\alpha+ \lambda)^{w(S_{1})} \lambda^{w(S_{2})}. 
$$

First, we prove $H_{\mathfrak{p}}^{(d)}(0) \neq 0$. 
Let us consider the constant term of $H_{\mathfrak{p}}^{(d)}(\lambda)$. 
If $S_{2} \neq \emptyset$, then $w(S_{2}) \neq 0$, and so the terms do not contribute to the constant term. 
When $S_{2}= \emptyset$, we get $w(S_{2})=0$, $L(S_{2}+1)= 1$ and $w(S_{1})= (q^{d}-1)/(q-1)$. %, that is, $w(S_{1})= q^{d-1}+ \cdots+ q+1$. 
Hence, $H_{\mathfrak{p}}^{(d)}(0)= \alpha^{(q^{d}-1)/(q-1)} \neq 0$. %, that is, $H_{\mathfrak{p}}^{(d)}(\lambda)$ does not have $0$ as root.  

\medskip

Next, we show that $\deg_{\lambda} H_{\mathfrak{p}}^{(d)}(\lambda)= (q^{d}-1)/(q-1)$. 
If $S_{2} \neq \emptyset$, then 
\begin{align*}
(q^{d}-1)/(q-1) &= q^{d-1}+ \cdots+ q+ 1 \\
&= w(S_{1})+ w(S_{2})+ w(S_{2}+1) \\
& > w(S_{1})+ w(S_{2}). 
\end{align*}
When $S_{2}= \emptyset$, we get $w(S_{1})= (q^{d}-1)/(q-1)$. 
Hence, $\deg_{\lambda} H_{\mathfrak{p}}^{(d)}(\lambda)= (q^{d}-1)/(q-1)$. 

\medskip

(c) \ We prove the separability of polynomials by using items (a) and (b). 
Recall that $g_{2d}= \lambda^{q^{2(d-1)}+ \cdots+ q^{2}+ 1}$ and $g_{2d+1}=0$. 
Suppose that $i= 2d+1$ for the recursion (\ref{recurrencegi}). 
Then, we can write 
\begin{align}\label{recurrenceg2d1}
g_{2d-1}^{q^{2}} \lambda+ g_{2d}^{q} \cdot (\alpha+ \lambda) = g_{2d-1} \lambda^{q^{2d-1}} + g_{2d} \cdot (\alpha+ \lambda)^{q^{2d}}. 
\end{align}
By differentiating both sides with respect to $\lambda$ and then multiplying both sides by $\lambda$, we have 
\begin{align}\label{recurrenceg2d1diff}
g_{2d-1}^{q^{2}} \lambda+ g_{2d}^{q} \lambda= g_{2d-1}^{\prime} \lambda^{q^{2d-1}+1} + g_{2d} \cdot (\alpha+ \lambda)^{q^{2d}}, 
\end{align}
where $g_{2d-1}^{\prime}$ is the derivative of $g_{2d-1}$. 
It follows from the equalities (\ref{recurrenceg2d1}) and (\ref{recurrenceg2d1diff}) that 
\begin{align}\label{recurrenceg2d1andd1diff}
g_{2d}^{q} \alpha= (g_{2d-1}- g_{2d-1}^{\prime} \lambda) \lambda^{q^{2d-1}}. 
\end{align}

% Recall that $g_{2d}= \lambda^{(q^{2d}-1)/(q^{2}-1)}$. 
Now, suppose that $H_{\mathfrak{p}}^{(d)}(\lambda)$ has a multiple root $\lambda_{0}$. 
Then, from item (a), $g_{2d-1}$ has the same multiple root $\lambda_{0}$, that is, $g_{2d-1}(\lambda_{0})= g_{2d-1}^{\prime}(\lambda_{0})=0$. 
Therefore, since the right-hand side of the equality (\ref{recurrenceg2d1andd1diff}) has the element $\lambda_{0}$ as a root, the left-hand side, too, has the element $\lambda_{0}$ as a root. 
Since the root of the left-hand side is $0$ only, we obtain $\lambda_{0}=0$. 
This contradicts item (b). 

\medskip

Next, suppose that $H_{\mathfrak{p}}^{(d)}(-\alpha^{q} s(s+1)^{q-1})$ has a multiple root $s_{0}$. 
Set $\lambda_{0}:= -\alpha^{q} s_{0}(s_{0}+1)^{q-1}$. 
Then
\begin{align*}
H_{\mathfrak{p}}^{(d)}(\lambda_{0}) &= 0, \qquad \text{and} \\
\frac{d}{ds} H_{\mathfrak{p}}^{(d)}(-\alpha^{q} s_{0}(s_{0}+1)^{q-1}) &= \frac{d}{d \lambda} H_{\mathfrak{p}}^{(d)}(\lambda_{0}) \cdot (-\alpha^{q} (s_{0}+1)^{q-2})=0.  
\end{align*}
Therefore, we have either $\frac{d}{d \lambda} H_{\mathfrak{p}}^{(d)}(\lambda_{0})=0$ or $-\alpha^{q} (s_{0}+1)^{q-2}=0$. 
If $\frac{d}{d \lambda} H_{\mathfrak{p}}^{(d)}(\lambda_{0})=0$, then $\lambda_{0}$ is a multiple root of $H_{\mathfrak{p}}^{(d)}(\lambda)$. 
However, $H_{\mathfrak{p}}^{(d)}(\lambda)$ does not have a multiple root. 
So we assume $-\alpha^{q} (s_{0}+1)^{q-2}=0$. %, that is, $s_{0}= -1$. 
Then, we obtain $\lambda_{0}= -\alpha^{q} s_{0}(s_{0}+1)^{q-1}=0$, which contradicts item (b). 
% Hence, $H_{\mathfrak{p}}^{(d)}(-\alpha^{q} s(s+1)^{q-1})$ is simple. 

\medskip

(d) \ 
Suppose that $\lambda_{0} \in \bar{\F}_{\mathfrak{p}}$ is any element such that $H_{\mathfrak{p}}^{(d)}(\lambda_{0}^{q+1})=0$, and let $\lambda:= \lambda_{0}^{q+1}$. 
Then, $\phi$ is supersingular at $\mathfrak{p}$ since $H_{\mathfrak{p}}^{(d)}(\lambda)=0$. 
It is known that a supersingular Drinfeld module is defined over $\F_{\mathfrak{p}}^{(2)}$, that is, $\lambda \in \F_{\mathfrak{p}}^{(2)}$ (see Remark 9.2 in \cite{Gek01}, Proposition 2.15 in \cite{Gek04}). 
Recall that $\alpha^{q^{d}}= \alpha$ and $\lambda^{q^{2d}}= \lambda$. 
Assume that $i=2d+1$ for the recursion (\ref{recurrencegi}). 
Then
$$
\lambda^{q(q^{2d}-1)/(q^{2}-1)} (\alpha+ \lambda)- \lambda^{(q^{2d}-1)/(q^{2}-1)} (\alpha+ \lambda)^{q^{2d}}= 0, 
$$
which can be written as, 
$$
\lambda^{(q^{2d}-1)/(q^{2}-1)} (\alpha+ \lambda)(\lambda^{(q^{2d}-1)/(q+1)}-1)= 0. 
$$
Hence, we have either $\lambda= -\alpha$ or $\lambda^{(q^{2d}-1)/(q+1)}= 1$. 
Notice that $\lambda \neq 0$ from item (b). 

First, consider the case where $\lambda^{(q^{2d}-1)/(q+1)}= 1$. 
With that, we get $\lambda_{0}^{q^{2d}-1}= \lambda^{(q^{2d}-1)/(q+1)}= 1$, and hence $\lambda_{0} \in \F_{\mathfrak{p}}^{(2)}$. 

Next, consider the case $\lambda= - \alpha$. 
When $d$ is odd, we get 
\begin{align*}
\lambda_{0}^{q^{2d}-1} &= \lambda_{0}^{(q^{d}+1)(q^{d}-1)}= (\lambda_{0}^{(q+1)(q^{d}-1)})^{(q^{d}+1)/(q+1)} \\
&= ((- \alpha)^{q^{d}-1})^{(q^{d}+1)/(q+1)}= 1, 
\end{align*}
and hence $\lambda_{0} \in \F_{\mathfrak{p}}^{(2)}$. 
If $d$ is even, then $\lambda= - \alpha$ is not a root of $H_{\mathfrak{p}}^{(d)}(\lambda)$. 
Indeed, we have 
\begin{align*}
H_{\mathfrak{p}}^{(d)}(-\alpha)= (-1)^{d/2} [1][3] \cdots [d-1] (-\alpha)^{q^{d-2}+ \cdots+ q^{2}+1} \neq 0
\end{align*} 
from Main theorem (1). 
\end{proof}
%%%%%%%%%

\medskip

%%%%%%%%%%
\begin{rem*}
(1) \ 
Much of the ideas of the above proof are similar to those in Section 5 in \cite{Gek83}. 
% That is, the type (\ref{drinfeldnormal}) different from the Gekeler's type (\ref{gekelernormal}) is merely used. 
The second claim of item (a) corresponds to the inseparability of Proposition 4.1 (c) in \cite{Gek91}. 

% The following fact was noticed by Gekeler: 
Gekeler noted the following on reading an earlier version of this paper (personal communication): 
The item (c) states that the supersingular locus on a moduli scheme is reduced, and this result has already been proven in a more general case.   

\medskip

(2) \ 
The background for item (b) is as follows. 
Let $\Sigma(\mathfrak{p})$ denote the set of supersingular points of $X(1)/\F_{\mathfrak{p}}$, that is, the set of supersingular $j$-invariants. 
It is known that the point $j= 0$ is supersingular if and only if $d$ is odd, and that 
$$
| \Sigma(\mathfrak{p}) |= 
\begin{cases}
\ \displaystyle \frac{q^{d}-1}{q^{2}-1} & \text{if $d$ is even}; \\
\ \displaystyle \frac{q^{d}-q}{q^{2}-1}+1 & \text{if $d$ is odd}
\end{cases}
$$
and $\Sigma(\mathfrak{p}) \subseteq \F_{\mathfrak{p}}^{(2)}$ (see Satz (5.9) in \cite{Gek83}, (2.14) in \cite{Gek04}). 
Notice that the covering $X_{0}(T) \to X(1)$ is given by $j= (\alpha+ \lambda)^{q+1}/\lambda$ (see the last part of Section \ref{secttower}). 
Since the point $j= \infty$ is not supersingular for $X(1)$, the point $\lambda=0$ is also not supersingular for $X_{0}(T)$. 
That is, $H_{\mathfrak{p}}^{(d)}(0) \neq 0$. 
Moreover, we can count the supersingular points of $X_{0}(T)/\F_{\mathfrak{p}}$ as follows: 

First, suppose that $d$ is even. 
Since all the supersingular points split completely in $X_{0}(T) \to X(1)$, the number is equal to 
$$
\frac{q^{d}-1}{q^{2}-1} \times (q+1)= \frac{q^{d}-1}{q-1}= \deg_{\lambda} H_{\mathfrak{p}}^{(d)}(\lambda). 
$$ 
Next, assume that $d$ is odd.
The supersingular point $j=0$ is totally ramified in $X_{0}(T) \to X(1)$ and the other supersingular points split completely in this covering. 
Hence, the number is equal to   
$$
\frac{q^{d}-q}{q^{2}-1} \times (q+1)+ 1= \frac{q^{d}-1}{q-1}= \deg_{\lambda} H_{\mathfrak{p}}^{(d)}(\lambda). 
$$

\medskip

(3) \ 
In Proposition 13 of \cite{BB}, Bassa and Beelen showed that a polynomial $p_{d}(s)$ is separable. 
Furthermore, in Theorem 18 of \cite{BB}, they proved that all roots of $p_{d}(s)$ are in $\F_{\mathfrak{p}}^{(2)}$. 
Although $p_{d}(s)$ is related to $H_{\mathfrak{p}}^{(d)}(\lambda)$ (see the remark in Section \ref{sectss}), their method of proof is completely different from the method used here. 
Our method of proof is straightforward. %very simple. 
\end{rem*}
%%%%%%%%%

%%%%%%%%%%%%%%%%%%%%%%%%%%%%%%%%%%%%%%%%%%
\section{An application for a supersingular polynomial}\label{secttower}
%%%%%%%%%%%%%%%%%%%%%%%%%%%%%%%%%%%%%%%%%%

In this section, we prove Corollary (Proposition \ref{maincor}), which was introduced in Section \ref{sectintro}. 
In the course of the proof, a polynomial identity (Proposition \ref{polyidentity} (b)) plays an important role. 
By combining this polynomial identity with a generalization of a result by Bezerra and Garcia in \cite{BG04}, the corollary is proven.   

\bigskip

% \medskip

Let $\mathfrak{p}= (p(T))$ denote a nonzero prime ideal of degree $d$ in $A= \F_{q}[T]$ such that $p(T) \neq T$. 
% Set $d= \deg p(T)$. 
Let $\alpha$ be a root of $p(T)$, and let $\F_{\mathfrak{p}}^{(2)}$ denote the quadratic extension of $\F_{\mathfrak{p}} := A/\mathfrak{p}= \F_{q}(\alpha)= \F_{q^{d}}$. 

Consider  a tower $\mathcal{E} = \mathcal{E}^{(d)} := (E_{0}, E_{1}, E_{2}, \ldots)$ that is recursively defined over $\F_{\mathfrak{p}}^{(2)}$ by the equation 
\begin{align}\label{BBequation}
Y(Y+1)^{q-1}= \frac{X^{q}}{(\alpha (X+1))^{q-1}} 
\end{align}
(see Definition 7.2.12 in \cite{Sti09}). 
This was first introduced by Bassa and Beelen in \cite{BB}. 

When $d=1$ (and so, without loss of generality, we can assume that $p(T)= T-1$ and $\alpha=1$). 
The tower $\mathcal{E}^{(1)}/\F_{q^{2}}$ was first introduced by Elkies (see the equation (25) in \cite{Elk01}), and was studied by Bezerra and Garcia (see the equation (1) in \cite{BG04}). 

Setting $X= (1-x)/x$ and $Y= (1-y)/y$ in the equation (\ref{BBequation}), we get 
\begin{align*} %\label{generalizedBG}
\frac{y-1}{y^{q}} = \frac{x^{q}-1}{x}, 
\end{align*}
which is the equation (1) in \cite{BG04}. 
Bezerra and Garcia showed that the genus of $\mathcal{E}^{(1)}/\F_{q^{2}}$ is 
\begin{align*}
\gamma(\mathcal{E}^{(1)}) := \lim_{n \to \infty} \frac{g(E_{n})}{q^{n}} = \frac{q}{q-1}, 
\end{align*}
where $g(E_{n})$ denotes the genus of $E_{n}$ (see Lemma 4 in \cite{BG04}). 

This result holds for any degree $d$, that is, $\gamma(\mathcal{E}^{(d)})= {q}/{(q-1)}$ (see the remark after Proposition \ref{maincor}, or Theorem 8.1 (iii) in \cite{Gek01}, Theorem 2.13 in \cite{Gek04}), and this is used in the proof of Corollary. 
% Because (???)

\bigskip

% \medskip

The following polynomial identities play an important role in our proof of Corollary. 

%%%%%%%%%%
\begin{prop}\label{polyidentity}
Assume that $A_{1}= \iota(T)+ \lambda$ and $A_{2}= \lambda$. 
Under this assumption, since the coefficients $\mathfrak{b}(d)$ are polynomials in $\lambda$, we can write $\mathfrak{b}(d)(\lambda):= \mathfrak{b}(d)$. 

\medskip

(a) \ 
Further, suppose that $L$ is an extension of $K$. 
Then 
\begin{align*}
& \mathfrak{b}{(d)}(-T^{q} s(s+1)^{q-1}) = (T(s+1))^{q^{d}-1} \mathfrak{b}{(d)} \left( \frac{-T^{q} s^{q}}{(T(s+1))^{q-1}} \right) \\
& + (T^{q^{d}-1}-1) (s+1)^{q^{d}-1} \mathfrak{b}{(d-1)} \left( \frac{-T^{q} s^{q}}{(T(s+1))^{q-1}} \right). 
\end{align*}

\medskip

(b) \ 
Instead, assume that $L$ is an extension of $\F_{\mathfrak{p}}$. 
Then 
\begin{align*}
H_{\mathfrak{p}}^{(d)}(-\alpha^{q} s(s+1)^{q-1}) &= (s+1)^{q^{d}-1} H_{\mathfrak{p}}^{(d)} \left( \frac{-\alpha^{q} s^{q}}{ (\alpha(s+1))^{q-1}} \right). %; \\
% H_{\mathfrak{p}}^{(d)}(-\alpha^{q} (s^{q}+s)^{q-1}) &= (s^{q-1}+1)^{q^{d}-1} H_{\mathfrak{p}}^{(d)} \left( -\alpha^{q} \left( \frac{s^{q}}{\alpha (s^{q-1}+1)} \right)^{q-1} \right). 
\end{align*}
\end{prop}
%%%%%%%%

\medskip

%%%%%%%%%%%
\begin{proof}
% Recall that $A/\mathfrak{p}= \F_{q}(\alpha)= \F_{q^{d}}$, and so $\alpha^{q^{d}-1}= 1$. 
% Recall that $\F_{q}(\alpha)= \F_{q^{d}}$ and $\alpha^{q^{d}-1}= 1$. 
Since $\alpha^{q^{d}-1}= 1$, item (b) follows from item (a). 
So, it is sufficient to prove item (a) only. %, which is similar to that of Theorem 15 in \cite{BB}. 
We show item (a) by induction on $d$. 
Notice that although item (b) is dependent on the prime ideal $\mathfrak{p}$, item (a) is not.  
% Furthermore, if we replace the variable $s$ in the first claim of the item (2) by $s^{q}$, then the second claim of the item (2) follows. 

\medskip

Recall that $\mathfrak{b}{(-1)}(\lambda)=0$, $\mathfrak{b}{(0)}(\lambda)=1$, and $\mathfrak{b}{(1)}(\lambda)= (-1) (\lambda/T^{q}+1)$ (see the definitions of $\mathfrak{a}(d)$ and $\mathfrak{b}(d)$ in Section \ref{sectss}). 
When $d=0$, the claim is obviously true. 
So, suppose that $d=1$. 
The left-hand side is then $s(s+1)^{q-1}-1$. 
The right-hand side is 
\begin{align*}
(T(s+1))^{q-1} \left( \frac{s^{q}}{(T(s+1))^{q-1}}-1 \right) &+ (T^{q-1}-1)(s+1)^{q-1} \\
& = s(s+1)^{q-1}-1. 
\end{align*}
Thus, the claim is true. 

For simplicity, let 
$$
L:= -T^{q}s(s+1)^{q-1} \qquad \text{and} \qquad R:= - \frac{T^{q} s^{q}}{(T(s+1))^{q-1}}. 
$$
Assume that the claim holds for $d-1$ and $d$. 
That is, assume that 
\begin{align*}
\mathfrak{b}{(d-1)}(L) &= (T(s+1))^{q^{d-1}-1} \mathfrak{b}{(d-1)}(R) \\
& + (T^{q^{d-1}-1}-1) (s+1)^{q^{d-1}-1} \mathfrak{b}{(d-2)}(R); \\
% \begin{align*}
% \end{align*}
\mathfrak{b}{(d)}(L) &= (T(s+1))^{q^{d}-1} \mathfrak{b}{(d)}(R) \\
& + (T^{q^{d}-1}-1) (s+1)^{q^{d}-1} \mathfrak{b}{(d-1)}(R). 
\end{align*}
It follows from the recursion (\ref{recursionb}) in Section \ref{sectss} that 
\begin{align}
\mathfrak{b}{(d+1)}(L) &= - \left( 1+ \left( \frac{L}{T^{q}} \right)^{q^{d}} \right) \mathfrak{b}{(d)}(L) \nonumber \\ 
& + (T^{1- q^{d}}-1) \left( \frac{L}{T^{q}} \right)^{q^{d}} \mathfrak{b}{(d-1)}(L);  \label{recursionbL} \\  
% \end{align}
% \begin{align}
\mathfrak{b}{(d)}(R) &= - \left( 1+ \left( \frac{R}{T^{q}} \right)^{q^{d-1}} \right) \mathfrak{b}{(d-1)}(R) \nonumber \\
& + (T^{1- q^{d-1}}-1) \left( \frac{R}{T^{q}} \right)^{q^{d-1}} \mathfrak{b}{(d-2)}(R); \label{recursionbR1} \\ 
% \end{align}
% \begin{align}
\mathfrak{b}{(d+1)}(R) &= - \left( 1+ \left( \frac{R}{T^{q}} \right)^{q^{d}} \right) \mathfrak{b}{(d)}(R) \nonumber \\
& + (T^{1- q^{d}}-1) \left( \frac{R}{T^{q}} \right)^{q^{d}} \mathfrak{b}{(d-1)}(R), \label{recursionbR2} 
\end{align} 
and these relations are used in this proof. 
First, we can compute the right-hand side of (\ref{recursionbL}) as follows: 
\begin{align*}
\mathfrak{b}{(d+1)}(L) &= - \left( 1+ \left( \frac{L}{T^{q}} \right)^{q^{d}} \right) \mathfrak{b}{(d)}(L) \\
& + (T^{1- q^{d}}-1) \left( \frac{L}{T^{q}} \right)^{q^{d}} \mathfrak{b}{(d-1)}(L) \\
&= \left( \left( s(s+1)^{q-1} \right)^{q^{d}}- 1 \right) \mathfrak{b}{(d)}(L) \\
& + (1- T^{1-q^{d}}) \left( s(s+1)^{q-1} \right)^{q^{d}} \mathfrak{b}{(d-1)}(L). 
\end{align*}
Then, it follows from the inductive hypothesis that 
\begin{align*}
\mathfrak{b}{(d+1)}(L) &= \left( \left( s(s+1)^{q-1} \right)^{q^{d}}- 1 \right) \mathfrak{b}{(d)}(L) \\
& + (1- T^{1-q^{d}}) \left( s(s+1)^{q-1} \right)^{q^{d}} \mathfrak{b}{(d-1)}(L) \\
&= A \cdot \mathfrak{b}{(d)}(R)+ B \cdot \mathfrak{b}{(d-1)}(R)+ C \cdot \mathfrak{b}{(d-2)}(R), 
\end{align*}
where 
\begin{align*}
A &:= \left( \left( s(s+1)^{q-1} \right)^{q^{d}}- 1 \right) (T(s+1))^{q^{d}-1}; \\
B &:= \left( \left( s(s+1)^{q-1} \right)^{q^{d}}- 1 \right) (T^{q^{d}-1}-1) (s+1)^{q^{d}-1} \\
& + (1- T^{1-q^{d}}) \left( s(s+1)^{q-1} \right)^{q^{d}}  (T(s+1))^{q^{d-1}-1}; \\
C &:= (1- T^{1-q^{d}}) \left( s(s+1)^{q-1} \right)^{q^{d}}  (T^{q^{d-1}-1}-1) (s+1)^{q^{d-1}-1}. 
\end{align*}
Next, by using  the equality (\ref{recursionbR1}), we calculate the right-hand side of the above equality and obtain 
\begin{align*}
& \mathfrak{b}{(d+1)}(L) = A \cdot \mathfrak{b}{(d)}(R)+ B \cdot \mathfrak{b}{(d-1)}(R)+ C \cdot \mathfrak{b}{(d-2)}(R) \\
&= \left( A+  (T^{q^{d}-1}-1) (s+1)^{q^{d+1}-1} \right) \mathfrak{b}{(d)}(R) \\
& + \left( B+ (T^{q^{d}-1}-1) (s+1)^{q^{d+1}-1} \left( 1+ \left( \frac{R}{T^{q}} \right)^{q^{d-1}} \right) \right) \mathfrak{b}{(d-1)}(R) \\
&= \left( A+  (T^{q^{d}-1}-1) (s+1)^{q^{d+1}-1} \right) \mathfrak{b}{(d)}(R) \\
& + (T^{q^{d}-1}-1) s^{q^{d+1}} (s+1)^{q^{d}-1}  \mathfrak{b}{(d-1)}(R). 
\end{align*}
Last, from the equality (\ref{recursionbR2}), we have 
\begin{align*}
& \mathfrak{b}{(d+1)}(L) =  \left( A+  (T^{q^{d}-1}-1) (s+1)^{q^{d+1}-1} \right) \mathfrak{b}{(d)}(R) \\
&+  (T^{q^{d}-1}-1) s^{q^{d+1}} (s+1)^{q^{d}-1}  \mathfrak{b}{(d-1)}(R) \\
&= (T(s+1))^{q^{d+1}-1}  \mathfrak{b}{(d+1)}(R) \\
&+ \left( A+  (T^{q^{d}-1}-1) (s+1)^{q^{d+1}-1}+ (T(s+1))^{q^{d+1}-1} \left( 1+ \left( \frac{R}{T^{q}} \right)^{q^{d}} \right)  \right) \mathfrak{b}{(d)}(R) \\
&=  (T(s+1))^{q^{d+1}-1}  \mathfrak{b}{(d+1)}(R)+  (T^{q^{d+1}-1}-1) (s+1)^{q^{d+1}-1} \mathfrak{b}{(d)}(R),  
\end{align*}
which is the desired result. 
\end{proof}
%%%%%%%%%

\medskip

%%%%%%%%%%
\begin{rem*}
% We give a remark about the above proposition. 
(1) \ 
The idea for the proof of item (a) is found in Theorem 15 in \cite{BB}. 

\medskip

(2) \ 
% In the above proposition, we want the item (2) rather than the item (1). 
In the course of the proof, the idea that the reduced polynomial $H_{\mathfrak{p}}^{(d)}(\lambda)$ is once-lifted to the unreduced polynomial $\mathfrak{b}{(d)}$ is very important. 
That is, item (b) can be proven by using item (a). 
In what follows, item (b) is often used, but item (a) is not used directly. 
% Because, I do not know a straightforward proof of the item (2). % the item (2) is proven by using the item (1). 
\end{rem*}
%%%%%%%%%

\medskip

As an application of Proposition \ref{polyidentity} (b), we have the following corollary. 

%%%%%%%%%
\begin{cor}\label{numberofss}
Let 
$$
\Omega= \Omega^{(d)} := \Big\{ s \in \bar{\F}_{\mathfrak{p}} \ \Big| \ H_{\mathfrak{p}}^{(d)}(-\alpha^{q} s(s+1)^{q-1})=0 \Big\}. 
$$
Then, $\Omega \subseteq \F_{\mathfrak{p}}^{(2)}$ and $|\Omega|= q \cdot \deg_{\lambda} H_{\mathfrak{p}}^{(d)}(\lambda)= q(q^{d}-1)/(q-1)$. 
\end{cor}
%%%%%%%%

\medskip

%%%%%%%%%%%
\begin{proof}
% Recall that $H_{\mathfrak{p}}^{(d)}(\lambda)$ does not have $\lambda= 0$ as root (see the proof of Proposition \ref{mainthm23} (2)). 
% It follows from Proposition \ref{mainthm23} (3) that if $H_{\mathfrak{p}}^{(d)}(\lambda)=0$, then $\lambda \in \F_{\mathfrak{p}}^{(2)}$. 
First, we show that $\F_{\mathfrak{p}}^{(2)}$ contains a primitive $(q+1)$th root of unity. 
If $\zeta \in \bar{\F}_{q}$ is a primitive $(q+1)$th root of unity, then $\zeta^{q^{2}-1}= (\zeta^{q+1})^{q-1}=1$, and so $\zeta \in \F_{q^{2}} \subseteq \F_{\mathfrak{p}}^{(2)}$. 

Next, we prove $-1 \notin \Omega$. 
If $-1 \in \Omega$, then $H_{\mathfrak{p}}^{(d)}(\lambda)$ has $\lambda= 0$ as a root, which contradicts Proposition \ref{mainthm23} (b). 
Hence, $-1 \notin \Omega$. 

For each $s \in \Omega$, we get 
\begin{align*}
(s+1)^{q^{d}-1} H_{\mathfrak{p}}^{(d)} \left( -\alpha^{q} \frac{s^{q}}{ (\alpha(s+1))^{q-1}} \right) &= H_{\mathfrak{p}}^{(d)}(-\alpha^{q} s(s+1)^{q-1}) \\
& = 0
\end{align*}
from Proposition \ref{polyidentity} (b). 
That is, 
\begin{align*}
H_{\mathfrak{p}}^{(d)}(-\alpha^{q} s(s+1)^{q-1}) &= 0; \\
H_{\mathfrak{p}}^{(d)} \left( -\alpha^{q} \frac{s^{q}}{ (\alpha(s+1))^{q-1}} \right) &= 0. 
\end{align*}
By Main theorem (2), we can write $\lambda^{q+1} = - \alpha^{q} s(s+1)^{q-1}$ and ${\lambda^{\prime}}^{q+1} = - \alpha^{q} s^{q}/(\alpha(s+1))^{q-1}$ for some $\lambda, \lambda^{\prime} \in \F_{\mathfrak{p}}^{(2)}$. 
Since $( \lambda \lambda^{\prime} )^{q+1}= \alpha^{q+1} s^{q+1}$ (that is, $s= \zeta \lambda \lambda^{\prime}/\alpha \in \F_{\mathfrak{p}}^{(2)}$), we have $\Omega \subseteq \F_{\mathfrak{p}}^{(2)}$. 

It follows from Proposition \ref{mainthm23} (c) that $H_{\mathfrak{p}}^{(d)}(-\alpha^{q} s(s+1)^{q-1})$ is separable. 
Thus, we obtain $|\Omega|= \deg_{s} (-\alpha^{q} s(s+1)^{q-1}) \cdot \deg_{\lambda} H_{\mathfrak{p}}^{(d)}(\lambda)= q(q^{d}-1)/(q-1)$ by Main theorem (2). 
\end{proof}
%%%%%%%%%

\medskip

%%%%%%%%%%
\begin{rem*}
(1) \ 
The set $\Omega$ can be identified with the set of supersingular points of the Drinfeld modular curve $X_{0}(T^{2})/\F_{\mathfrak{p}}$. 
This can be proven as follows. 

Recall that the set $S$ of roots of $H_{\mathfrak{p}}^{(d)}(\lambda)$ is the set of supersingular points of $X_{0}(T)/\F_{\mathfrak{p}}$. 
It is known that the covering $X_{0}(T^{2}) \to X_{0}(T)$ is defined by $\lambda= -\alpha^{q} s (s+1)^{q-1}$ (see the last part of this section). 
Hence, $\Omega$ is the set of supersingular points of $X_{0}(T^{2})$. 
Since all the supersingular points of $X_{0}(T)$ split completely in this covering, and $S \subseteq \F_{\mathfrak{p}}^{(2)}$ (see Proposition \ref{mainthm23} (d)), we have $\Omega \subseteq \F_{\mathfrak{p}}^{(2)}$. 

\medskip

(2) \ 
In Corollary 19 in \cite{BB}, Bassa and Beelen proved the same result as the above corollary for another polynomial $p_{d}(s(s+1)^{q-1})$. 
Their method of proof is completely different from the method used here. 
% However, their method of proof is completely different. 
\end{rem*}
%%%%%%%%%

\bigskip

% \medskip

The following proposition follows as a corollary (which we call Corollary). %  introduced in Section \ref{sectintro}. 
Here is a justification. 
Although Corollary is expressed in terms of the curves $X_{0}(T^{n})$, the proposition is expressed in terms of the function fields $E_{n}$. 
It is well-known that the nonsingular complete curves can be put in one-to-one correspondence with the function fields of one variable (see, for example, Corollary 6.12 in \cite{Har77}, Remark 2.5 of Chapter II in \cite{Sil09}).  
Under this correspondence, a function field of the curve $X_{0}(T^{n+2})$ corresponds exactly to the function field $E_{n}$ (see the last part of this section). % given by the equation (\ref{BBequation}) 

Hence, it is sufficient to prove the following proposition. 

%%%%%%%%%%
\begin{prop}\label{maincor}
The tower $\mathcal{E}/\F_{\mathfrak{p}}^{(2)}$ is asymptotically optimal. 
That is, $\lambda(\mathcal{E})= q^{d}-1$. 
\end{prop}
%%%%%%%%

\medskip

%%%%%%%%%%
\begin{proof}
Recall that $\gamma(\mathcal{E}) = {q}/{(q-1)}$ (cf. Lemma 4 in \cite{BG04}). 

Now, we compute the limit $\nu(\mathcal{E}) := \lim_{n \to \infty} N(E_{n}/\F_{\mathfrak{p}}^{(2)})/q^{n}$. 
We write the zero of $x_{0}-a$ in the rational function field $E_{0}= \F_{\mathfrak{p}}^{(2)}(x_{0})$ as $P_{a}:= P_{x_{0}-a}$. 

First, we show that for each $a \in \Omega$, the place $P_{a}$ splits completely in $\mathcal{E}/\F_{\mathfrak{p}}^{(2)}$. 
Recall that $-1 \notin \Omega$. 
Fix an element $b \in \bar{\F}_{\mathfrak{p}}$ such that $b(b+1)^{q-1}= a^{q}/(\alpha (a+1))^{q-1}$. 
It follows from Proposition \ref{polyidentity} (b) that 
\begin{align*}
H_{\mathfrak{p}}^{(d)}(-\alpha^{q} b(b+1)^{q-1}) &= H_{\mathfrak{p}}^{(d)} \left( -\alpha^{q} \frac{a^{q}}{ (\alpha(a+1))^{q-1}} \right) \\
& = H_{\mathfrak{p}}^{(d)}(-\alpha^{q} a(a+1)^{q-1})/(a+1)^{q^{d}-1} \\
&= 0, 
\end{align*}
and so $b \in \Omega$. 
Thus, we get $b \in \F_{\mathfrak{p}}^{(2)}$ by Corollary \ref{numberofss}. 
By Kummer's theorem (see Theorem 3.3.7 in \cite{Sti09}),  the place $P_{a}$ splits completely in $E_{1}/E_{0}$. 
The desired assertion follows by induction. 

Next, we calculate $\nu(\mathcal{E})$. 
By the first claim, we have $N(E_{n}/\F_{\mathfrak{p}}^{(2)}) \geq |\Omega| \cdot q^{n}$, and so 
\begin{align*}
\nu(\mathcal{E}) &= \lim_{n \to \infty} \frac{N(E_{n}/\F_{\mathfrak{p}}^{(2)})}{q^{n}} \\
& \geq |\Omega| = \frac{q(q^{d}-1)}{q-1}
\end{align*}
from Corollary \ref{numberofss}. 

Combining the above results, we get $\lambda(\mathcal{E})= \nu(\mathcal{E})/\gamma(\mathcal{E}) \geq q^{d}-1$. 
It then follows from the Drinfeld-Vl\u{a}du\c{t} bound (see Theorem 7.1.3 in \cite{Sti09}) that $q^{d}- 1 \leq \lambda(\mathcal{E}) \leq A(q^{2d})= q^{d}-1$. 
\end{proof}
%%%%%%%

\medskip

%%%%%%%%%%
\begin{rem*}
(1) \ 
Recall that $d= \deg_{T} p(T)$ is the degree of $\mathfrak{p}= (p(T))$. 
Our corollary is a generalization of Theorem 1 of \cite{BG04} (the Bezerra-Garcia theorem), generalizing that result to allow an arbitrary degree $d$. 
That is, the case when $d=1$ corresponds exactly to that theorem. 
Although our corollary is a special case of Theorem 2.16 in \cite{Gek04} and Theorem 4.2.38 in \cite{TV91}, our proof is more elementary, and explicitly describes the set $\Omega$ of degree-one places. 
For this reason, our result has applications to coding theory (see \cite{Sti06}, Chapters 7 and 8 in \cite{Sti09}, Parts 3 and 4 in \cite{TV91}, Chapters 3 and 4 in \cite{TVD07}). 
The Bezerra-Garcia theorem is a special case of Theorem 10.1 of Gekeler in \cite{Gek01}. 

\medskip

First, we compare Corollary \ref{numberofss} with the result by Bezerra and Garcia (and a result of Elkies in \cite{Elk01}) in terms of the splitting locus $\text{Split}(\mathcal{E})$ (see Definition 7.2.9 (a) in \cite{Sti09}). 
% A huge difference between Theorem 2.16 in \cite{Gek04} and Theorem 10.1 in \cite{Gek01} is the splitting rate. Assume that $d=1$. 
Recall that when $d=1$, we have $p(T)=T-1$. 
Bezerra and Garcia showed that the places corresponding to the roots of $x^{q}+ x-1=0$ split completely (see Page 152 in \cite{BG04}). 
Setting $s= (1-x)/x$, we get 
\begin{align*}
x^{q}+ x- 1= 0 \quad & \Leftrightarrow \quad -s(s+1)^{q-1}+1= 0 \\
& \Leftrightarrow \quad H_{\mathfrak{p}}^{(1)}(-s(s+1)^{q-1})= 0, 
\end{align*}
and so the special case ($d=1$) of Corollary \ref{numberofss} coincides with the Bezerra-Garcia result. 
Furthermore, the set of roots for the second equation coincides with the set (26) given by Elkies in \cite{Elk01}. 

Now, we explain a difference between Corollary \ref{numberofss} and the Bezerra-Garcia result. 
Recall that $\Sigma(\mathfrak{p})$ denotes the set of supersingular points of $X(1)/\F_{\mathfrak{p}}$. 
When $d=1$, we know that $| \Sigma(\mathfrak{p}) |= 1$ and that the point $j= 0$ is the supersingular point of $X(1)$. % (and is also elliptic point). 
Since the covering $X_{0}(T^{2}) \to X_{0}(T)$ is given by $\lambda= - \alpha^{q} s (s+1)^{q-1}$ and the covering $X_{0}(T) \to X(1)$ is defined by $j= (\alpha+ \lambda)^{q+1}/\lambda$ (see the last part of this section), we easily see that all the roots of $-s(s+1)^{q-1}+1= 0$ are above $j= 0$. 
% Considering more cases, 
When $d \geq 3$ is odd, we know that $| \Sigma(\mathfrak{p}) | \geq q+1$. 
Therefore, the completely splitting points are above several supersingular points. 
In fact, the points $j= 0$ and $j= [1] (1- (\alpha^{q}- \beta)^{q-1})$ ($\beta \in \F_{q}$) are supersingular points (see Proposition 16 in \cite{Sch95}). 
Hence, the splitting locus $\text{Split}(\mathcal{E}^{(d)})$ cannot be computed by using the approach employed in \cite{BG04}. 

\medskip

Next, we compare our result with the result by Bezerra and Garcia in terms of the ramification locus $\text{Ram}(\mathcal{E})$ (see Definition 7.2.9 (b) in \cite{Sti09}). 
%Suppose that $d=1$. 
Bezerra and Garcia showed that the places corresponding to $x=0$, $x=1$ or $x= \infty$ are ramified, and the other places are unramified (see Section 3 in \cite{BG04}). 
Setting $s= (1-x)/x$, we get 
\begin{align*}
x=0, \ x=1, \ x= \infty \quad & \Leftrightarrow \quad s= \infty, \ s= 0, \ s= -1, 
\end{align*}
respectively. 
% It is known that the point $j= \infty$ is the cusp of $X(1)$. 
% It is known that the covering $X_{0}(T^{2}) \to X_{0}(T)$ is given by $\lambda= -\alpha^{q} s (s+1)^{q-1}$, and that the covering $X_{0}(T) \to X(1)$ is defined by $j= (\alpha+ \lambda)^{q+1}/\lambda$. 
We can see that the points $s=0,-1,\infty$ are above $j= \infty$. 
Similarly, in our case, we know that the ramified points are above $j= \infty$. 
Hence, the ramification locus $\text{Ram}(\mathcal{E}^{(d)})$ can be calculated using the method of \cite{BG04}.  
% Notice that the ramified points $s= \infty, 0, -1$ are independent of the constant field $\F_{\mathfrak{p}}^{(2)}$. 

\medskip

(2) \ 
In the proof of Proposition \ref{maincor}, as a consequence, we see that 
\begin{align*}
\nu(\mathcal{E})= \lim_{n \to \infty} \frac{N(E_{n}/\F_{\mathfrak{p}}^{(2)})}{q^{n}} = |\Omega|, 
\end{align*}
and so $N(E_{n}/\F_{\mathfrak{p}}^{(2)})= |\Omega|+ o(q^{n})$ ($n \to \infty$). 
\end{rem*}
%%%%%%%%%

% \bigskip
% Let $\mathcal{F}:= \mathcal{F}_{\alpha}= (F_{i})$ be a tower over $\F_{\mathfrak{p}}^{(2)}$ which is recursively defined by the equation 
% \begin{align}\label{GSequation} y^{q}+ y= \frac{x^{q}}{\alpha (x^{q-1}+1)}, \end{align}
% which was first introduced by Bassa and Beelen in \cite{BB}. 
% When $p(T)= T-1$, the tower $\mathcal{F}_{1}/\F_{q^{2}}$ was studied by Garcia and Stichtenoth in \cite{GS96}. 
% The tower $\mathcal{F}$ is a supertower of our tower $\mathcal{E}$. 
% Similarly, $\mathcal{E}$ is one of subtower of $\mathcal{F}$. 
% This is proven as follows: 
% Set $X= x^{q-1}$ and $Y= y^{q-1}$ in the equation (\ref{BBequation}). 
% Then, we get  
% \begin{align*} y^{q-1}(y^{q-1}+1)^{q-1}= \frac{x^{q(q-1)}}{(\alpha (x^{q-1}+1))^{q-1}}, \end{align*}
% that is, 
% \begin{align*} (y^{q}+ y)^{q-1}= \left( \frac{x^{q}}{\alpha (x^{q-1}+1)} \right)^{q-1}, \end{align*}
% which is the equation (\ref{GSequation}). 

% \bigskip
% Bassa-Beelen (arXiv) and Schweizer (95) ($(X_{0}(T^{n+2}))_{n \geq 0}$?)

\bigskip

% \medskip

Finally, we consider a background of the tower $\mathcal{E}$ in terms of the Drinfeld modular curves $X_{0}(T^{n})$. 

First, we consider the genus, which yields another proof that $\gamma(\mathcal{E})= q/(q-1)$. 
It follows from Theorem 8.1 (iii) in \cite{Gek01} (or Theorem 2.13 in \cite{Gek04}) that the genus of $X_{0}(T^{n})$ is given by 
\begin{align*}
g(X_{0}(T^{n})) &= \frac{q^{n-1}- q^{\lceil (n-1)/2 \rceil}- q^{\lceil (n-2)/2 \rceil}+ 1}{q-1} \\
&= \frac{q}{q-1}+ o(q^{n-2}) \qquad (n \to \infty). 
\end{align*}
Hence, $g(X(1))= g(X_{0}(T))= g(X_{0}(T^{2}))=0$ and $g(X_{0}(T^{n})) \geq 1$ for $n \geq 3$. 
This is used below. 

\medskip

Next, we consider the origin of the equation (\ref{BBequation}), which is due to Sections 2 and 3 in \cite{BB}. 
Let $\phi$ be a rank-$2$ Drinfeld module over $K$, defined by $\phi_{T}= T+ (T+ \lambda_{0}) \tau+ \lambda_{0} \tau^{2}$ with $j$-invariant $j_{0}= (T+ \lambda_{0})^{q+1}/\lambda_{0}$. 
It is known that its $T$-isogenous (rank-$2$) Drinfeld module $\phi^{\prime}$ is given by $\phi_{T}^{\prime}= T+ (T^{q}+ \lambda_{0}) \tau+ \lambda_{0}^{q} \tau^{2}$, and its $j$-invariant is $j_{1}(T)= (T^{q}+ \lambda_{0})^{q+1}/\lambda_{0}^{q}$. 
% Assume that $j_{0}$ is transcendental over $K$. 
% Notice that $j_{0}$ is independent of $T$, but $j_{1}(T)$ depends on $T$. 
There is some modular polynomial $\Phi_{T}(X,Y) \in A[X,Y]$ such that $\Phi_{T}(j_{0}, j_{1}(T))= 0$, and this polynomial is very complicated (see any of \cite{Sch95}, \cite{BL97}, \cite{CHJ06}, \cite{BB11, BB12}). 
% First, we construct the function field of $X_{0}(T^{n})$. 
% For a transcendental element $j_{0}$, the equation $\Phi_{T}(j_{0}, j_{1}(T))=0$ gives a plane (singular) model for the Drinfeld modular curve $X_{0}(T)$. 

\medskip

Second basement: 
Assume that $j_{0}$ is transcendental over $K$. 
That is, assume that the function field $K_{0}:= K(j_{0})$ is rational, which is a function field of $X(1)$. 
The Drinfeld modular curve $X_{0}(T)$ is defined by the equation $\Phi_{T}(j_{0}, j_{1}(T))= 0$, and its function field $K(X_{0}(T))$ is given as $K_{1}:= K(j_{0}, j_{1}(T))$. 
Recall that the relation $\Phi_{T}(j_{0}, j_{1}(T))= 0$ provides a plane model for $X_{0}(T)$, and that 
$$
[K_{1}: K_{0}]= q^{\deg (T)} \prod_{P \mid T, \atop \text{$P$ is prime}} \left( 1+ \frac{1}{q^{\deg (P)}} \right)= q+1. 
$$
It is known that, for each integer $n \geq 0$, the function field $K(X_{0}(T^{n}))$ can be written as % of the curve $X_{0}(T^{n})$ is given as 
$$
K_{n}:= K(j_{0}, j_{1}(T), j_{1}(T^{2}), \ldots, j_{1}(T^{n})), 
$$
where 
$$
\Phi_{T}(j_{0}, j_{1}(T))= 0 \qquad \text{and} \qquad \Phi_{T}(j_{1}(T^{i}), j_{1}(T^{i+1}))=0
$$
for $1 \leq i<n$. 
Notice that $\Phi_{T}(j_{1}(T^{i}), Y)$ is reducible over $K_{i}$, and that $[K_{n}: K_{0}]= (q+1) q^{n-1}$. 

\medskip

First basement: 
Since $X_{0}(T)$ is also of genus $0$, its function field $K(X_{0}(T))$ is rational. 
In fact, $K(X_{0}(T))$ can be given as $K_{1}= K(\lambda_{0})$ by using the transcendental element $\lambda_{0}$ (see Proposition 3 in \cite{Sch95}). 
The function field $K(X_{0}(T^{2}))$ is then given as $K_{2}= K(\lambda_{0}, \lambda_{1})$, where 
$$
\frac{(T+ \lambda_{1})^{q+1}}{\lambda_{1}}= \frac{(T^{q}+ \lambda_{0})^{q+1}}{\lambda_{0}^{q}}. % \ \left( = j_{1}(T) \right). 
$$
% In other words, the element $u_{1}$ is a generator of $K(j_{1}(T), j_{1}(T^{2}))$. 
This relation is not minimal with respect to degree, in the following sense: 
Since 
\begin{align*}
& \frac{(T+ \lambda_{1})^{q+1}}{\lambda_{1}}- \frac{(T^{q}+ \lambda_{0})^{q+1}}{\lambda_{0}^{q}} \\
& = \frac{T^{q+1}- \lambda_{0} \lambda_{1}}{\lambda_{0}^{q} \lambda_{1}} \left( T^{q^{2}}+ \lambda_{0}^{q}- \left( T^{q+1}- \lambda_{0} \lambda_{1} \right)^{q-1} \left( T+ \lambda_{1} \right) \right) 
\end{align*}
and $\lambda_{0} \lambda_{1} \neq T^{q+1}$, we can obtain a new relation
\begin{align*}
T^{q^{2}}+ \lambda_{0}^{q}- \left( T^{q+1}- \lambda_{0} \lambda_{1} \right)^{q-1} \left( T+ \lambda_{1} \right)= 0. 
\end{align*}
That is, 
\begin{align*}
\frac{T^{q}+ \lambda_{0}}{T+ \lambda_{1}}= \left( \frac{T^{q+1}- \lambda_{0} \lambda_{1}}{T^{q}+ \lambda_{0}} \right)^{q-1}, 
\end{align*}
which is a minimal relation (with respect to degree). % between $\lambda_{0}$ and $\lambda_{1}$. 
This minimal relation is used below. 

For each integer $n \geq 1$, the function field $K(X_{0}(T^{n}))$ is 
\begin{align*}
K_{n}= K(\lambda_{0}, \lambda_{1}, \ldots, \lambda_{n-1}), 
\end{align*}
where 
\begin{align*}
\frac{(T+ \lambda_{i+1})^{q+1}}{\lambda_{i+1}}= \frac{(T^{q}+ \lambda_{i})^{q+1}}{\lambda_{i}^{q}}
\end{align*}
for $0 \leq i<n-1$. 
Notice that $[K_{1}: K_{0}]= q+1$ and $[K_{n}: K_{1}]= q^{n}$. 

\medskip

Ground floor: 
Since $X_{0}(T^{2})$ is also of genus $0$, its function field $K(X_{0}(T^{2}))$ is rational. 
In fact, $K(X_{0}(T^{2}))$ can be written as $K_{2}= K(s_{0})$ by using the transcendental element $s_{0}$ such that  
$$
-T \cdot (s_{0}+1) = \frac{T^{q+1}- \lambda_{0} \lambda_{1}}{T^{q}+ \lambda_{0}}
$$
(see Section 3 in \cite{BB}). 
Since 
\begin{align*}
\lambda_{0}= - T^{q} s_{0} (s_{0}+ 1)^{q-1} \qquad \text{and} \qquad \lambda_{1}= - \frac{T^{q} s_{0}^{q}}{(T (s_{0}+1))^{q-1}}, 
\end{align*}
we obtain $K_{1}= K(\lambda_{0}, \lambda_{1})= K(s_{0})$. 

For each integer $n \geq 2$, the function field $K(X_{0}(T^{n}))$ is 
\begin{align*}
K_{n}= K(s_{0}, s_{1}, \ldots, s_{n-2}), 
\end{align*}
where 
\begin{align*}
s_{i+1} (s_{i+1}+ 1)^{q-1}= \frac{s_{i}^{q}}{(T (s_{i}+1))^{q-1}} 
\end{align*}
for $0 \leq i <n-2$, which is just the equation (\ref{BBequation}). 

By using the technique of Elkies in \cite{Elk98, Elk01}, the sequence of $X_{0}(T^{n})/\F_{\mathfrak{p}}^{(2)}$ ($n \geq 2$) corresponds to the tower $\mathcal{E}$. 
In particular, a function field of $X_{0}(T^{n+2})/\F_{\mathfrak{p}}^{(2)}$ corresponds to the function field $E_{n}$. 

%%%%%%%%%%%%%%%%%%%%%%%
\begin{acknowledgements*}
The author is grateful to Professor Ernst-Ulrich Gekeler, who pointed out several mistaken citations in Section 1 after the author uploaded an earlier version of this paper to arXiv. 
This work was supported by JSPS KAKENHI Grant Number 15K17508. 
% The author thank the referee for his or her valuable comments and careful review of this paper. 
\end{acknowledgements*}
%%%%%%%%%%%%%%%%%%%%%

% \bigskip Hypergeometric function (???)

%%%%%%%%%%%%%%%%%%%%%

%%%%%%%%%%%%%%%%%

\end{document}